\documentclass[11 pt,english]{article}
\usepackage[activeacute]{babel}
\usepackage{amsfonts} 
\usepackage{mathrsfs} 
\usepackage{amssymb, amsmath, amsfonts}
\usepackage{makeidx}
\usepackage{epsfig}
\usepackage{color}
\usepackage[T1]{fontenc}
\usepackage{graphics,graphicx,graphpap}
\usepackage[ansinew]{inputenc}
\usepackage{amsthm}
\usepackage{subfigure}
\usepackage{tikz}

\usetikzlibrary{calc, positioning, through, math, arrows, backgrounds}
\tikzstyle{vtx}=[inner sep=1pt,draw, shape=circle, font=\tiny]
\tikzstyle{line}=[inner sep=3pt,draw, shape=rectangle, line width = 3pt]
\tikzset{>=stealth}
\tikzstyle{lbl}=[inner sep = 1 pt, fill = white, midway]

\newcommand\blfootnote[1]{%
  \begingroup
  \renewcommand\thefootnote{}\footnote{#1}%
  \addtocounter{footnote}{-1}%
  \endgroup
}

\if@twoside \oddsidemargin 6pt \evensidemargin 6pt \marginparwidth 20pt
 \else \oddsidemargin 6pt \evensidemargin 6pt \marginparwidth 20pt \fi
\newcommand{\ZZ}{\mathbb{Z}}

\newcommand{\G}{\Gamma}

\newcommand{\cP}{\mathcal{P}}
\newcommand{\cR}{\mathcal{R}}

\newcommand{\Wr}{\mathop{\rm W}}
\newcommand{\QWr}{\mathop{\rm QW}}

\newtheorem{theorem}{Theorem}[section]
\newtheorem{proposition}[theorem]{Proposition}
\newtheorem{corollary}[theorem]{Corollary}
\newtheorem{lemma}[theorem]{Lemma}
\newtheorem{question}[theorem]{Question}

\theoremstyle{definition}

\newtheorem{problem}[theorem]{Problem}
\newtheorem{construction}[theorem]{Construction}

\begin{document}

\begin{center}
\Large{\textbf{Self-reverse labelings of distance magic graphs}} \\ [+4ex]
\end{center}

\begin{center}
Petr Kov\'a\v{r}{\small $^{a}$},\   
Ksenija Rozman{\small $^{b,c}$},\
Primo\v z \v Sparl{\small $^{b,d,e,*}$}
\\

\medskip
{\it {\small
$^a$Technical University of Ostrava, Department of Applied Mathematics, Ostrava, Czech Republic\\
$^b$Institute of Mathematics, Physics and Mechanics, Ljubljana, Slovenia\\
$^c$University of Primorska, FAMNIT, Koper, Slovenia\\
$^d$University of Ljubljana, Faculty of Education, Ljubljana, Slovenia\\
$^e$University of Primorska, Institute Andrej Maru\v si\v c, Koper, Slovenia\\
}}
\end{center}

\blfootnote{
Email addresses: 
petr.kovar@vsb.cz (P.~Kov\'a\v{r}), ksenija.rozman@pef.uni-lj.si (K.~Rozman), primoz.sparl@pef.uni-lj.si (P. \v Sparl)
\\
* - corresponding author
}


\hrule

\begin{abstract}
A graph is distance magic if it admits a bijective labeling of its vertices by integers from $1$ up to the order of the graph in such a way that the sum of the labels of all the neighbors of a vertex is independent of a given vertex. We introduce the concept of a self-reverse distance magic labeling of a regular graph which allows for a more compact description of the graph and the labeling in terms of the corresponding quotient graph. We show that the members of several known infinite families of tetravalent distance magic graphs admit such labelings. We present a novel general construction producing a new distance magic graph from two existing ones. Using it we show that for each integer $n \geq 6$, except for the odd integers up to $19$, there exists a connected tetravalent graph of order $n$ admitting a self-reverse distance magic labeling. We also determine all connected tetravalent graphs up to order $30$ admitting a self-reverse distance magic labeling. The obtained data suggests a number of natural interesting questions giving several possibilities for future research.
\end{abstract}
\hrule

\begin{quotation}
\end{quotation}

\section{Introduction}
\label{sec:Intro}

A graph $\G$ of order $n$ is {\em distance magic} if it admits a bijective labeling of its vertices with integers from $1$ up to $n$ such that the sum of the labels of the neighbors of a vertex is independent of the given vertex. Following~\cite{SugFroMilRyaWal09} such a labeling is nowadays called a distance magic labeling of $\G$. It was pointed out in~\cite{MikSpa21} however, that when working with regular graphs a somewhat different definition of a distance magic labeling can be used, leading to the same definition of a regular graph being distance magic. In this alternative definition one lets the labels be the integers of the arithmetic progression from $1-n$ to $n-1$ with common difference $2$ and insists that the sum of the labels of the neighbors of each vertex equals $0$ (see Section~\ref{sec:Prelim} for details). Since this alternative definition is often more convenient to work with, we will be using it throughout this paper. 

Ever since Vilfred~\cite{Vil94} introduced the concept of distance magic graphs (using a somewhat different terminology --- the term distance magic was first used in~\cite{SugFroMilRyaWal09}) researchers have been studying this intriguing class of graphs from various points of view, as is testified by the number of papers on the topic (see the dynamic survey on graph labeling maintained by Gallian~\cite{GalD}). While there do exist some papers considering general properties of distance magic graphs (such as for instance~\cite{AruKamVij14, ONeSla27}, where it is proved that the so-called magic constant of a distance magic graph is unique, or~\cite{SugFroMilRyaWal09}, in which some basic necessary conditions for a graph to be distance magic are given), most of the papers focus on classifying distance magic members of some well-known basic families of graphs, or on finding sufficient or necessary conditions for some of the most standard products of graphs (Cartesian, direct, lexicographic or strong) to be distance magic (see~\cite{AnhCicPetTep15, CicFro16, CicFroKroRad16, GreKov13, MikSpa21, MikSpa24, TiaHouHouGao21} for some examples, but refer to \cite{GalD} for the long list of the corresponding references). The lack of attempts at a systematic investigation and/or development of general methods for the study of distance magic graphs is most probably due to the fact that the whole class of distance magic graphs is quite rich and extremely diverse and thus does not seem to allow for a universal approach. As a consequence, many researchers have started to investigate various related concepts in recent years, such as for instance closed distance magic graphs and group distance magic graphs (see~\cite{GalD} for a long list of references on the topic). 

Despite what was just said we are convinced that the topic of distance magic graphs provides interesting possibilities for general investigations, especially if one restricts to regular graphs and even more so if one further restricts to those with a large degree of symmetry. That this is indeed the case is indicated by a recent result from~\cite{MikSpa21}, which shows that for regular graphs the property of being distance magic can be stated in terms of eigenvalues and eigenvectors of the corresponding adjacency matrix. This new insight enabled several recent results (see for instance~\cite{FerMalMikRaz23, MikSpa21,MikSpa23, MikSpa24, RozSpa24}).

In the present paper we propose another such direction of research by introducing the concept of a self-reverse distance magic labeling of a regular graph (the concept is briefly explained at the end of this section, but see Section~\ref{sec:Self} for the definition). While there exist regular distance magic graphs admitting no self-reverse distance magic labeling, it is interesting to note that for most of the infinite families of regular graphs whose members are known to be distance magic, the graphs in fact admit self-reverse distance magic labelings. Moreover, there exist infinite families of regular distance magic graphs for which all distance magic labelings are self-reverse (see Sections~\ref{sec:tetravalent} and~\ref{sec:VT} where this is discussed in detail). 

While the concept of a self-reverse distance magic labeling seems to be interesting in its own (we point out a number of natural questions that arise in the subsequent sections), one of the main benefits of studying these labelings is that they allow for a somewhat more compact description of the graph and the labeling in question, and at the same time link the investigation of distance magic graphs to the well-studied theory of graph covers~\cite{GroTuc77}. In this way one can hope to develop a nice tool for constructing distance magic regular graphs from smaller graphs. 
\medskip

One of our results is a complete classification of all orders for which a tetravalent distance magic graph admitting a self-reverse distance magic labeling exists (see Theorem~\ref{the:classification}). Another result of the research leading to this manuscript is a complete list of all self-reverse distance magic labelings of connected tetravalent graphs of orders up to $30$. The obtained data suggests a number of very natural and interesting questions (see~Section~\ref{sec:VT}), therefore giving several possibilities for future research. In particular, tetravalent vertex-transitive graphs (see Section~\ref{sec:VT} for a definition) admitting a (self-reverse) distance magic labeling do exist, but seem to be extremely rare. This calls for a thorough analysis of tetravalent distance magic graphs admitting a large degree of symmetry in the future. As a curiosity, we mention that among the small tetravalent distance magic graphs admitting a self-reverse distance magic labeling, very interesting corresponding quotient graphs are encountered. For instance, one of the two such ``nontrivial'' tetravalent graphs of order $20$ (see Section~\ref{sec:classification}) yields the Petersen graph (with a semiedge at each vertex) as the corresponding quotient graph (see Figure~\ref{fig:VTexamples}).

We conclude this introductory section with a brief description of the main idea behind self-reverse distance magic labelings of regular graphs and the outline of the paper. Suppose that $\G$ is a distance magic graph of order $n$ and that $\ell$ is a corresponding distance magic labeling (in the classical sense or the one from Section~\ref{sec:Prelim} which we use in this paper). Then $\ell$ induces a natural partition $\cP_\ell$ of the vertex set of $\G$ into pairs (and one singleton in the case that $n$ is odd) in which the vertices $u$ and $v$ are in the same set from $\cP_\ell$ if and only if for some $i$ they have the $i$-th smallest and the $i$-th largest label. One can then consider the corresponding quotient graph $\G_\ell$ of $\G$ with respect to $\cP_\ell$ in which two vertices are adjacent whenever there exists at least one edge in $\G$ having its two endvertices in the two corresponding partition sets. In general much of the information on $\G$ (and the labeling $\ell$) is lost by considering just $\G_\ell$. If however the labeling $\ell$ is particularly nice in the sense that between any two ``adjacent'' sets from $\cP_\ell$ we either have just a perfect matching or we have all possible edges (in which case $\ell$ is said to be {\em self-reverse} --- see Section~\ref{sec:Self} for a precise definition), then the situation is significantly different. In this case the graph $\G$ is what is known in the literature as a regular $2$-fold cover (at least in the case that $n$ is even) over the corresponding natural generalization of the above mentioned quotient graph, and so all of the information on $\G$ can be encoded via this quotient graph and the corresponding voltages. Moreover, the labeling $\ell$ itself can be described uniquely via a labeling of the quotient graph and these voltages. As is explained in Section~\ref{sec:Self}, in this case exchanging each label by its ``reverse'' (the ``antipodal label'') yields an equivalent distance magic labeling of the graph, thereby explaining why we say it is self-reverse. 

The paper is structured as follows. After first gathering the basic notions, making conventions and fixing the notation that we use throughout the paper, we introduce the concept of a self-reverse distance magic labeling of a regular graph in Section~\ref{sec:Self}. In Section~\ref{sec:tetravalent} we then focus on tetravalent graphs. We introduce the concept of degeneracy and show that the so-called wreath graphs are the only examples of tetravalent graphs admitting degenerate self-reverse distance magic labelings. We also determine the ``types'' of self-reverse distance magic labelings of wreath graphs and introduce the concept of the quotient with respect to a self-reverse distance magic labeling. In Section~\ref{sec:construction} we introduce a novel general construction yielding a new regular distance magic graph from two existing ones (see Proposition~\ref{pro:con}) and show that under some additional assumptions  the resulting graph admits a self-reverse distance magic labeling (Corollary~\ref{cor:SRcon}). Together with obtained computational data this enables us to classify all orders for which a connected tetravalent graph admitting a self-reverse distance magic labeling exists (see Theorem~\ref{the:classification}). In Section~\ref{sec:VT} we point out a number of interesting natural questions suggested by the obtained data. This opens possibilities for further investigations on (self-reverse) distance magic labelings of regular graphs.

\section{Preliminaries}
\label{sec:Prelim}

In this section we gather the notation, definitions and conventions that we will be using throughout the paper. 

All graphs considered in this paper will be finite, connected and undirected. Moreover, except in the case of quotient graphs, where we will also allow semiedges, all graphs will be simple. For a graph $\G$ with vertex set $V(\G)$, we denote the fact that the vertices $u, v \in V(\G)$ are adjacent in $\G$ by $u \sim_\G v$ (or simply by $u \sim v$ if the graph $\G$ is clear from the context). We denote the corresponding edge by $uv$. We let $\G(v)$ be the neighbourhood of $v$ and let $\G[v]$ be the closed neighbourhood of $v$ (that is, $\G[v] = \G(v) \cup \{v\}$).

In this paper we focus on regular graphs. As mentioned in Section~\ref{sec:Intro}, a nonstandard definition of a distance magic labeling is convenient to work with in such a case. Let $\G$ be a regular graph of order $n$ and let $\mathcal{I}_n = \{1-n, 3-n, \ldots , n-1\}$ be the set consisting of the $n$ elements forming the arithmetic progression from $1-n$ to $n-1$ with common difference $2$. Then a bijective mapping $\ell \colon V(\G) \to \mathcal{I}_n$ is a {\em distance magic labeling} of $\G$, whenever for each vertex the sum of the labels of its neighbors equals $0$. Clearly, if $\ell \colon V(\G) \to \mathcal{I}_n$ is some mapping then the mapping $\ell' \colon V(\G) \to \{1,2,\ldots , n\}$, defined by $\ell'(v) = (1+n+\ell(v))/2$ for all $v \in V(\G)$, is a distance magic labeling in the ``classical sense'' if and only if $\ell$ is a distance magic labeling (according to the above definition). Therefore, the above alternative definition does indeed lead to an equivalent definition of a regular graph being distance magic. We point out one more time that throughout the rest of this paper we will constantly be working with this alternative definition, so if we ever want to refer to the definition from~\cite{SugFroMilRyaWal09} (where the set of the labels is $\{1,2,\ldots , n\}$), we will speak of a ``distance magic labeling in the classical sense''. 

Since the so-called wreath graphs, which constitute one of the most natural families of tetravalent distance magic graphs, will play an important role in our investigation of self-reverse distance magic labelings, we review the definition for self-completeness. For an integer $n \geq 3$, the {\em wreath graph} $\Wr(n)$ is the tetravalent graph of order $2n$ with vertex set $\{x_i \colon i \in \ZZ_n\} \cup \{y_i \colon i \in \ZZ_n\}$ (where $\ZZ_n$ denotes the ring of residue classes modulo $n$), in which for each $i \in \ZZ_n$ each of $x_i$ and $y_i$ is adjacent to each of $x_{i-1}, y_{i-1}, x_{i+1}, y_{i+1}$ (all subscripts are to be computed modulo $n$). We mention that the graph $\Wr(n)$ can of course be viewed as the wreath product (also called the lexicographic product) of the cycle of order $n$ by the edgeless graph of order $2$, indicating why it is called this way.

\section{The concept of a self-reverse distance magic labeling}
\label{sec:Self}

Let us now introduce the concept of a self-reverse distance magic labeling of a regular graph. Let $\G$ be a regular distance magic graph, let $n$ be its order and let $\ell$ be a corresponding distance magic labeling, where the labels are the elements of $\mathcal{I}_n$. For each vertex $v$ of $\G$ let $v^\ell$ be the unique vertex of $\G$ such that $\ell(v) + \ell(v^\ell) = 0$. Observe that if $n$ is odd then there exists precisely one vertex $v$ such that $v^{\ell} = v$ (the one with $\ell(v) = 0$) --- we call this $v$ the {\em central vertex} of $\G$ {\em with respect to} $\ell$ in this case. We let $\cP_\ell$ be the partition $\cP_\ell = \{\{v, v^{\ell}\} \colon v \in V(\G)\}$ of $V(\G)$ and call it the {\em partition corresponding to} $\ell$. 

The distance magic labeling $\ell$ of $\G$ is said to be {\em self-reverse} if for each pair of vertices $u, v$ of $\G$ the vertices $u$ and $v$ are adjacent if and only if the vertices $u^\ell$ and $v^\ell$ are adjacent. In other words, $\ell$ is self-reverse if and only if the involutory permutation of the vertex set of $\G$ interchanging each vertex $v$ with its paired vertex $v^\ell$ is an automorphism of $\G$. Clearly, this can be rephrased in yet another way, as is recorded in the following straightforward observation.

\begin{lemma}
\label{le:obser}
Let $\G$ be a regular distance magic graph, let $\ell$ be a corresponding distance magic labeling and let $\cP_\ell$ be the partition of $V(\G)$ corresponding to $\ell$. Then the following three statements are equivalent. 
\begin{itemize}
\itemsep = 0pt
	\item[(i)] The labeling $\ell$ is self-reverse.
	\item[(ii)] The permutation of $V(\G)$, interchanging each $v \in V(\G)$ with $v^\ell$, is an automorphism of $\G$.
	\item[(iii)] For each pair of distinct sets $A, B \in \cP_\ell$, the bipartite subgraph of $\G$ with vertex set $A \cup B$ and consisting of all the edges of $\G$ with one endvertex in $A$ and the other in $B$, has no edges, or is a complete bipartite graph, or consists of a pair of disjoint edges.
\end{itemize}
\end{lemma}

Observe that the natural distance magic labeling of the wreath graph $\Wr(m)$, $m \geq 3$, where one sets $\ell(x_i) = 2i+1$ and $\ell(y_i) = -(2i+1)$ for each $i \in \{0,1,\ldots, m-1\}$, is self-reverse, showing that such labelings indeed do exist. In fact, we shall prove (see Proposition~\ref{pro:wreath}) that whenever $m$ is not divisible by $4$, each distance magic labeling of the wreath graph $\Wr(m)$ is self-reverse. On the other hand, it should come as no surprise that there do exist regular distance magic graphs having no self-reverse distance magic labeling and that there also exist regular distance magic graphs admitting a self-reverse and also a non-self-reverse distance magic labeling (see Section~\ref{sec:tetravalent}). In fact, for most of the infinite families of tetravalent graphs which are known to be distance magic, its members admit self-reverse distance magic labelings (see Section~\ref{sec:VT} where this is discussed in some detail).

Before explaining our choice for the name self-reverse we remark that in~\cite{AnhCicPetTep15} a somewhat related concept was studied. There the authors studied so-called balanced (distance magic) labelings of regular graphs of even order. In our terminology a balanced distance magic labeling of a regular graph of even order is a self-reverse distance magic labeling $\ell$ such that for any two distinct sets $A, B \in \cP_\ell$ the corresponding bipartite subgraph from the above lemma either has no edges or is a complete bipartite graph. As we show in Section~\ref{sec:tetravalent} (see Proposition~\ref{pro:degenerate}), at least in the case of tetravalent graphs, all such examples are in fact wreath graphs and will thus not be of great interest to us.

\begin{lemma}
\label{le:reverse}
Let $\G$ be a regular distance magic graph and let $\ell$ be a corresponding distance magic labeling of $\G$. Then setting $\ell'(v) = -\ell(v)$ for each $v \in V(\G)$, we obtain another distance magic labeling of $\G$.
\end{lemma}

Now, let $\G$ be a regular distance magic graph, let $\ell$ be a corresponding distance magic labeling and let $\ell'$ be as in Lemma~\ref{le:reverse} (we call it the {\em reverse} of $\ell$). By definition the labelings $\ell$ and $\ell'$ are of course two different (distance magic) labelings of $\G$. However, it seems reasonable to say that two labelings of a graph $\G$ are {\em equivalent} (or essentially the same) whenever for each pair of labels $i$ and $j$ the vertices labeled $i$ and $j$ by $\ell$ are adjacent if and only if the vertices labeled $i$ and $j$ by $\ell'$ are also adjacent. The following proposition, whose easy proof is left to the reader, finally explains why a self-reverse distance magic labeling of a regular graph is called this way.

\begin{proposition}
\label{pro:SRname}
Let $\G$ be a regular distance magic graph, let $\ell$ be a corresponding distance magic labeling of $\G$ and let $\ell'$ be its reverse from Lemma~\ref{le:reverse}. Then $\ell'$ is equivalent to $\ell$ if and only if $\ell$ is self-reverse.
\end{proposition}

\section{Self-reverse distance magic labelings of tetravalent graphs}
\label{sec:tetravalent}

While one can of course study self-reverse labelings of regular graphs in full generality, the focus of this paper is on tetravalent graphs. First and foremost because this is the smallest interesting valence for regular distance magic graphs, but also since for tetravalent graphs the self-reverse labelings come in two essentially different types, as we show next. 

Before stating the result we introduce the following terminology. Suppose $\ell$ is a distance magic labeling of a tetravalent graph $\G$. Then $\ell$ is said to be {\em degenerate} if there exist vertices $u$ and $v$ such that $u \neq u^\ell$, $v \neq v^\ell$ (which is always the case if $\G$ is of even order) and $u$ is adjacent to both $v$ and $v^\ell$ (note that this automatically implies that $v \notin \{u, u^\ell\}$). In other words, $\ell$ is degenerate if there exist distinct elements of the corresponding partition $\cP_\ell$, none of which is a singleton, such that there exist two incident edges between them.

\begin{proposition}
\label{pro:degenerate}
Let $\G$ be a connected tetravalent distance magic graph admitting a self-reverse distance magic labeling $\ell$. If $\ell$ is degenerate, then for each vertex $v$ of $\G$ the vertices $v$ and $v^\ell$ are distinct and have the same neighborhoods. Consequently, $\G$ is of even order and is a wreath graph.
\end{proposition}

\begin{proof}
Suppose $\ell$ is degenerate and let $u$ and $v$ be adjacent vertices such that $u \neq u^\ell$, $v \neq v^\ell$ and $u \sim  v^\ell$. Let $w$ and $w'$ be the remaining neighbors of $u$, besides $v$ and $v^\ell$. Since $\ell$ is a distance magic labeling and $\ell(v) + \ell(v^\ell) = 0$, it follows that $\ell(w) + \ell(w') = 0$ as well, and so $w' = w^\ell$. In particular, $w^\ell \neq w$. Since $\ell$ is self-reverse, $u^\ell$ is adjacent to each of $v, v^\ell, w$ and $w^\ell$, and so $u$ and $u^\ell$ have the same neighborhoods. Since $\G$ is connected, an inductive argument shows that for each vertex $x$ of $\G$ the vertices $x$ and $x^\ell$ are distinct and have the same neighborhood. It is now clear that $\G$ is a wreath graph and is thus of even order.
\end{proof}

The above proposition shows that, at least in the case of tetravalent graphs, the only balanced (in the sense of~\cite{AnhCicPetTep15}) connected tetravalent distance magic graphs are the wreath graphs. The interesting question is thus whether there are any tetravalent distance magic graphs admitting a non-degenerate self-reverse distance magic labeling. As we will see, this is indeed the case. In fact, for some well-known ``non-wreath'' tetravalent distance magic graphs the only distance magic labelings are self-reverse (see Section~\ref{sec:VT}). But before we switch our attention to ``non-wreath'' graphs let us first briefly revisit the wreath graphs to see what kind of distance magic labelings they admit.

\begin{proposition}
\label{pro:wreath}
Let $m \geq 3$ be an integer. If $m$ is not divisible by $4$, then each distance magic labeling of the wreath graph $\Wr(m)$ is self-reverse and is degenerate. For each $m \geq 4$ with $4 \mid m$ the graph $\Wr(m)$ admits a degenerate self-reverse distance magic labeling and a non-degenerate self-reverse distance magic labeling. In addition, if $m \geq 8$, it also admits a non-self-reverse distance magic labeling.
\end{proposition}

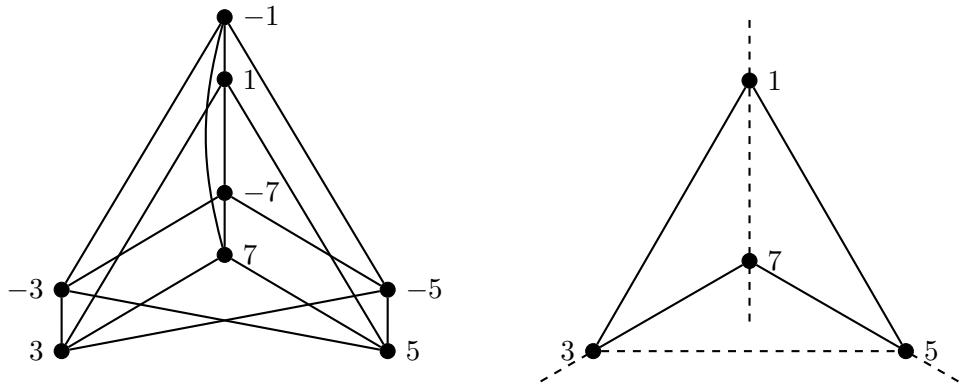
\begin{figure}[htbp]
\begin{center}
\subfigure{
\begin{tikzpicture}[scale = .55]
\node[vtx, fill=black, inner sep = 2pt, label=0:{$-1$}] (A0) at (90:5) {};
\node[vtx, fill=black, inner sep = 2pt, label=0:{$1$}] (B0) at (90:3.5) {};
\node[vtx, fill=black, inner sep = 2pt, label=180:{$-3$}] (A1) at (202:4.25) {};
\node[vtx, fill=black, inner sep = 2pt, label=180:{$3$}] (B1) at (218:5) {};
\node[vtx, fill=black, inner sep = 2pt, label=0:{$-5$}] (A2) at (338:4.25) {};
\node[vtx, fill=black, inner sep = 2pt, label=0:{$5$}] (B2) at (322:5) {};
\node[vtx, fill=black, inner sep = 2pt, label=0:{$-7$}] (A3) at (90:0.75) {};
\node[vtx, fill=black, inner sep = 2pt, label=0:{$7$}] (B3) at (270:0.75) {};
\draw[thick] (A0) -- (A1);
\draw[thick] (B0) -- (B1);
\draw[thick] (A0) -- (A2);
\draw[thick] (B0) -- (B2);
\draw[thick] (A1) -- (A3);
\draw[thick] (B1) -- (B3);
\draw[thick] (A2) -- (A3);
\draw[thick] (B2) -- (B3);
\draw[thick] (A1) -- (B2);
\draw[thick] (B1) -- (A2);
\draw[thick] (A0) -- (B0);
\draw[thick] (A1) -- (B1);
\draw[thick] (A2) -- (B2);
\draw[thick] (A3) -- (B3);
\draw[thick] (B0) to (A3);
\draw[thick] (A0) to [bend right = 15] (B3);
\end{tikzpicture}
}
\hspace{8mm}
\subfigure{
\begin{tikzpicture}[scale = .6]
\foreach \j in {0,1,2}{
\pgfmathtruncatemacro{\k}{2*\j+1}
\node[vtx, fill=black, inner sep = 2pt, label=180*\j:{$\k$}] (A\j) at (360*\j/3 + 90: 4) {};
}
\node[vtx, fill=black, inner sep = 2pt, label=0:{$7$}] (A3) at (0:0) {};
\draw[thick] (A0) -- (A1);
\draw[thick] (A0) -- (A2);
\draw[dashed, thick] (A0) -- (A3);
\draw[dashed, thick] (A1) -- (A2);
\draw[thick] (A1) -- (A3);
\draw[thick] (A2) -- (A3);

\begin{scope}[on background layer]
\draw[dashed, thick, shorten <= -7mm] (A0.90) -- (A0.270);
\draw[dashed, thick, shorten <= -7mm] (A1.210) -- (A1.30);
\draw[dashed, thick, shorten <= -7mm] (A2.330) -- (A2.150);
\draw[dashed, thick, shorten <= -7mm] (A3.270) -- (A3.90);
\end{scope}
\end{tikzpicture}
}
\caption{A non-degenerate self-reverse labeling of $\Wr(4)$ and the corresponding quotient.}
\label{fig:W(4)quo}
\end{center}
\end{figure}

\begin{proof}
Let $\ell$ be a distance magic labeling of $\Wr(m)$. For each $i \in \ZZ_m$ let $\ell_i = \ell(x_i) + \ell(y_i)$. Considering the neighbors of $x_{i+1}$ we find that $\ell_{i+2} = -\ell_{i}$ for all $i \in \ZZ_m$ (computation of indices is performed modulo $m$). In particular, 
\begin{equation}
\label{eq:wreaths}
	\ell_{i+4} = \ell_i\ \text{for all}\ i \in \ZZ_m.
\end{equation} 
Of course, $\sum_{i\in \ZZ_m}\ell_i = 0$, and so if $m$ is odd, $\ell_i = 0$ holds for all $i \in \ZZ_m$. Similarly, if $m \equiv 2 \pmod{4}$, then~\eqref{eq:wreaths} implies that $\ell_{i+2} = \ell_i$ for all $i \in \ZZ_m$, and then the above observation that $\ell_{i+2} = -\ell_i$ again implies that $\ell_i = 0$ for all $i \in \ZZ_m$. Therefore, $\ell(x_i) + \ell(y_i) = 0$ for all $i \in \ZZ_m$, and so the labeling $\ell$ is clearly self-reverse and degenerate.

The graph $\Wr(4)$ is isomorphic to the complete bipartite graph $K_{4,4}$. It is easy to see that for a distance magic labeling $\ell$ of this graph either the vertices labeled $i$ and $-i$ are contained in the same bipartition set for all $i \in \{1,3,5,7\}$, or this holds for no such $i$. It is now clear that in any case $\ell$ is self-reverse and that in the former case it is degenerate, while in the latter case it is not (see the left part of Figure~\ref{fig:W(4)quo}).

Finally, suppose $m = 4m_0$ for some $m_0 \geq 2$. A degenerate self-reverse distance magic labeling of $\Wr(m)$ is obtained by setting $\{\ell(x_i), \ell(y_i)\} = \{2i-2m+1, 2m-2i-1\}$ for all $i \in \{0,1,\ldots , m-1\}$. A non-degenerate self-reverse distance magic labeling of $\Wr(m)$ can be obtained as follows. For each $i$ with $0 \leq i < m_0$, set 
$$
	\begin{array}{ccl}
	\{\ell(x_{2i}), \ell(y_{2i})\} & = & \{(-1)^i (8i+1), (-1)^{i+1}(8i+3)\},\\
	\{\ell(x_{-1-2i}), \ell(y_{-1-2i})\} & = & \{(-1)^{i+1} (8i+1), (-1)^{i}(8i+3)\},\\
	\{\ell(x_{2i+1}), \ell(y_{2i+1})\} & = & \{(-1)^i (8i+5), (-1)^{i+1}(8i+7)\},\\
	\{\ell(x_{-2-2i}), \ell(y_{-2-2i})\} & = & \{(-1)^{i+1} (8i+5), (-1)^{i}(8i+7)\}.
	\end{array}
$$
We leave it to the reader to verify that this is indeed a distance magic labeling of $\Wr(m)$ and that it is self-reverse but is not degenerate. Letting $\tilde{\ell}$ be the labeling obtained from the above labeling $\ell$ by only changing the labels of $x_0, y_0, x_1$ and $y_1$ so that $\{\tilde{\ell}(x_0), \tilde{\ell}(y_0)\} = \{5,-7\}$ and $\{\tilde{\ell}(x_1), \tilde{\ell}(y_1)\} = \{1,-3\}$, we clearly get a new distance magic labeling of $\Wr(m)$, which however is not self-reverse since $x_1 \sim x_2$ but $x_1^{\tilde{\ell}}$ (which is one of $x_{-1}$ and $y_{-1}$) is not adjacent to $x_2^{\tilde{\ell}}$ (which is one of $x_{-3}$ and $y_{-3}$).
\end{proof}

We now explain how for a tetravalent graph $\G$ admitting a non-degenerate self-reverse distance magic labeling $\ell$, the graph $\G$ and the labeling $\ell$ can be described in a concise way. Let $\cP_\ell$ be the partition of the vertex set of $\G$ corresponding to the labeling $\ell$. 

Suppose first that $\G$ is of even order. In this case the elements of $\cP_\ell$ are unordered pairs of distinct vertices, and we can represent $\G$ and $\ell$ as follows. We construct a labeled graph (with possible semiedges) having one vertex for each element of $\cP_\ell$ and where the vertex corresponding to $\{v, v^\ell\}$ is labeled by $|\ell(v)|$. If $v \sim v^\ell$, we put a (dashed) semiedge at the vertex representing $\{v, v^\ell\}$. Moreover, for any pair of adjacent vertices $u$ and $v$ of $\G$ such that $v \neq u^\ell$ we put a solid or dashed edge between the vertices corresponding to $\{u,u^\ell\}$ and $\{v, v^\ell\}$, depending on whether $\ell(u)$ and $\ell(v)$ have the same sign or not, respectively. We call the constructed labeled graph with colored edges (solid and dashed) the {\em quotient} of $\G$ corresponding to $\ell$. In Figure~\ref{fig:W(4)quo}, the non-degenerate self-reverse distance magic labeling for the graph $\Wr(4)$ from the proof of Proposition~\ref{pro:wreath} and the corresponding quotient are presented. We point out that in the language of graph covers~\cite{GroTuc77, MalNedSko00}, the graph $\G$ is a regular $\ZZ_2$-cover over this quotient where the solid edges correspond to the edges with voltage $0$ while the dashed ones correspond to those with voltage $1$.

\begin{figure}[hbtp]
\begin{center}
\subfigure{
\begin{tikzpicture}[scale = .55]
\node[vtx, fill=black, inner sep = 2pt, label=0:{$0$}] (A0) at (90:4) {};
\node[vtx, fill=black, inner sep = 2pt, label=162:{$2$}] (A1) at (162:5) {};
\node[vtx, fill=black, inner sep = 2pt, outer sep = 0pt, label=90:{$-2$}] (B1) at (162:3.5) {};
\node[vtx, fill=black, inner sep = 2pt, label=198:{$10$}] (A2) at (198:5) {};
\node[vtx, fill=black, inner sep = 2pt, outer sep = -2pt, label=170:{$-10$}] (B2) at (198:3.5) {};
\node[vtx, fill=black, inner sep = 2pt, label=234:{$14$}] (A3) at (234:8) {};
\node[vtx, fill=black, inner sep = 2pt, label=54:{$-14$}] (B3) at (234:6.5) {};
\node[vtx, fill=black, inner sep = 2pt, outer sep = -2pt, label=20:{$6$}] (C3) at (234:4) {};
\node[vtx, fill=black, inner sep = 2pt, label=90:{$-6$}] (D3) at (234:2.5) {};
\node[vtx, fill=black, inner sep = 2pt, label=270:{$20$}] (A4) at (270:8) {};
\node[vtx, fill=black, inner sep = 2pt, label=90:{$-20$}] (B4) at (270:6.5) {};
\node[vtx, fill=black, inner sep = 2pt, label=270:{$16$}] (C4) at (270:4) {};
\node[vtx, fill=black, inner sep = 2pt, label=90:{$-16$}] (D4) at (270:2.5) {};
\node[vtx, fill=black, inner sep = 2pt, label=306:{$12$}] (A5) at (306:8) {};
\node[vtx, fill=black, inner sep = 2pt, label=126:{$-12$}] (B5) at (306:6.5) {};
\node[vtx, fill=black, inner sep = 2pt, label=270:{$4$}] (C5) at (306:4) {};
\node[vtx, fill=black, inner sep = 2pt, outer sep = -2pt, label=90:{$-4$}] (D5) at (306:2.5) {};
\node[vtx, fill=black, inner sep = 2pt, label=342:{$8$}] (A6) at (342:5) {};
\node[vtx, fill=black, inner sep = 2pt, outer sep = -1pt, label=100:{$-8$}] (B6) at (342:3.5) {};
\node[vtx, fill=black, inner sep = 2pt, label=18:{$18$}] (A7) at (18:5) {};
\node[vtx, fill=black, inner sep = 2pt, label=198:{$-18$}] (B7) at (18:3.5) {};

\draw[thick] (A0) -- (A1);
\draw[thick] (A0) -- (B1);
\draw[thick] (A0) -- (A7);
\draw[thick] (A0) -- (B7);

\draw[thick] (A1) -- (A2);
\draw[thick] (B1) -- (B2);
\draw[thick] (A1) to [bend right = 12] (A6);
\draw[thick] (B1) -- (B6);
\draw[thick] (A1) -- (B7);
\draw[thick] (B1) -- (A7);

\draw[thick] (A2) to [bend right = 12] (A7);
\draw[thick] (B2) -- (B7);
\draw[thick] (A2) -- (B3);
\draw[thick] (B2) -- (A3);
\draw[thick] (A2) -- (D3);
\draw[thick] (B2) -- (C3);

\draw[thick] (A3) -- (A4);
\draw[thick] (B3) -- (B4);
\draw[thick] (A3) to [bend left = 25] (C3);
\draw[thick] (B3) to [bend left = 25] (D3);
\draw[thick] (C3) -- (C4);
\draw[thick] (D3) -- (D4);
\draw[thick] (A3) -- (D4);
\draw[thick] (B3) -- (C4);
\draw[thick] (C3) -- (B4);
\draw[thick] (D3) -- (A4);

\draw[thick] (A4) -- (C5);
\draw[thick] (B4) -- (D5);
\draw[thick] (C4) -- (A5);
\draw[thick] (D4) -- (B5);
\draw[thick] (C4) -- (D5);
\draw[thick] (D4) -- (C5);
\draw[thick] (A4) -- (B5);
\draw[thick] (B4) -- (A5);
\draw[thick] (A5) to [bend right = 25] (D5);
\draw[thick] (B5) to [bend right = 25] (C5);

\draw[thick] (A5) -- (A6);
\draw[thick] (B5) -- (B6);
\draw[thick] (C5) -- (A6);
\draw[thick] (D5) -- (B6);

\draw[thick] (A6) -- (B7);
\draw[thick] (B6) -- (A7);
\end{tikzpicture}
}
\hspace{3mm}
\subfigure{
\begin{tikzpicture}[scale = .55]
\node[vtx, fill=black, inner sep = 2pt, label=0:{$0$}] (A0) at (90:5) {};
\node[vtx, fill=black, inner sep = 2pt, label=162:{$2$}] (A1) at (162:5) {};
\node[vtx, fill=black, inner sep = 2pt, label=198:{$10$}] (A2) at (198:5) {};
\node[vtx, fill=black, inner sep = 2pt, label=234:{$14$}] (A3) at (234:6) {};
\node[vtx, fill=black, inner sep = 2pt, label=54:{$6$}] (C3) at (234:3) {};
\node[vtx, fill=black, inner sep = 2pt, label=270:{$20$}] (A4) at (270:6) {};
\node[vtx, fill=black, inner sep = 2pt, label=90:{$16$}] (C4) at (270:3) {};
\node[vtx, fill=black, inner sep = 2pt, label=306:{$12$}] (A5) at (306:6) {};
\node[vtx, fill=black, inner sep = 2pt, label=126:{$4$}] (C5) at (306:3) {};
\node[vtx, fill=black, inner sep = 2pt, label=342:{$8$}] (A6) at (342:5) {};
\node[vtx, fill=black, inner sep = 2pt, label=18:{$18$}] (A7) at (18:5) {};
\draw[thick] (A0) -- (A1);
\draw[thick] (A0) -- (A7);
\draw[dashed, thick] (A1) -- (A7);
\draw[thick] (A1) -- (A6);
\draw[thick] (A1) -- (A2);
\draw[thick] (A2) -- (A7);
\draw[dashed, thick] (A2) -- (A3);
\draw[dashed, thick] (A2) -- (C3);
\draw[thick] (A3) -- (C3);
\draw[thick] (A3) -- (A4);
\draw[thick] (C3) -- (C4);
\draw[dashed, thick] (A3) -- (C4);
\draw[dashed, thick] (C3) -- (A4);
\draw[dashed, thick] (A4) -- (A5);
\draw[dashed, thick] (C4) -- (C5);
\draw[thick] (A4) -- (C5);
\draw[thick] (C4) -- (A5);
\draw[dashed, thick] (A5) -- (C5);
\draw[thick] (A6) -- (C5);
\draw[thick] (A6) -- (A5);
\draw[dashed, thick] (A6) -- (A7);
\end{tikzpicture}
}
\caption{A self-reverse labeling of a distance magic tetravalent graph of order $21$ and the corresponding quotient.}
\label{fig:Odd21quo}
\end{center}
\end{figure}
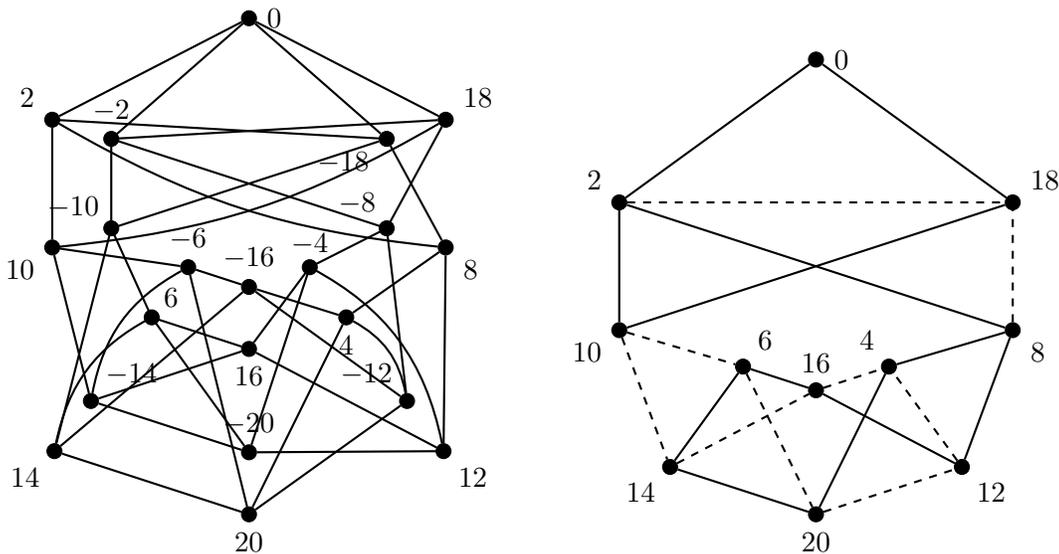

When the order of $\G$ is odd, the situation is slightly different, since in this case the central vertex with respect to $\ell$ (the vertex $v$ with $\ell(v) = 0$) gives rise to a singleton in $\cP_\ell$. The corresponding quotient is then constructed in the same way as in the case of graphs of even order, with the following exception. The central vertex $v$ is adjacent to some $u$ and $u^\ell$, and some $w$ and $w^\ell$ with $u \notin \{w, w^\ell\}$. In the quotient graph we simply connect the vertex corresponding to the singleton $\{v\}$ from $\cP_\ell$ to each of the vertices corresponding to $\{u,u^\ell\}$ and $\{w,w^\ell\}$, to each with a single solid edge. Since the vertex corresponding to $\{v\}$ will be labeled $0$, the fact that this vertex is ``special'' can be reconstructed from its label as well as from the fact that it is the only vertex in the quotient having only two neighbors. An example of a self-reverse distance magic labeling and a corresponding quotient for a distance magic tetravalent graph of order $21$ is given in Figure~\ref{fig:Odd21quo}.

We conclude this section with the following remarks, providing examples of connected tetravalent graphs that admit no self-reverse distance magic labelings. In~\cite{RozSpa25}, connected distance magic tetravalent graphs of small orders were analyzed. It was established that up to order $16$ only three examples which are not wreath graphs exist. Two of them are the so-called quasi-wreath graphs $\mathrm{QW}(3,3)$ (depicted on the right in Figure~\ref{fig:merged_example}) and $\mathrm{QW}(7)$ of orders $12$ and $14$, respectively (see~\cite{RozSpa25} for the definition and details), while the third one is a graph of order $16$, obtained from $\mathrm{QW}(7)$ using a simple construction from~\cite{KovFroKov12}. The distance magic members of the family of quasi-wreath graphs were classified and it was proved that there are infinitely many of them. However, it can be deduced from the proof of~\cite[Proposition~4.3]{RozSpa25}, that no distance magic quasi-wreath graph (except for the smallest one $\mathrm{QW}(3)$, which is isomorphic to $\mathrm{W}(3)$) admits a self-reverse distance magic labeling (with the notation from~\cite{RozSpa25}, one simply has to look at the vertices $x_2$ and $y_2$ and their neighbors). 

This thus confirms that we indeed have all possibilities: there are infinitely many connected tetravalent distance magic graphs that do not admit a self-reverse distance magic labeling, there are infinitely many connected tetravalent distance magic graphs for which all distance magic labelings are self-reverse, and there are also infinitely many connected tetravalent distance magic graphs, that admit self-reverse distance magic labelings, as well as non-self-reverse distance magic labelings.

\section{A general construction}
\label{sec:construction}

In this section we present a general construction that, given two distance magic graphs satisfying certain nice enough properties, yields a larger distance magic graph obtained from these two graphs. We think this construction is interesting in its own and deserves to be further explored, but in this paper it will mainly be used to help us classify all orders for which a connected tetravalent distance magic graph admitting a (non-degenerate) self-reverse distance magic labeling exists (see the next section).
Before presenting the construction we introduce some terminology and notation that we will be using. 

Let $\G$ be a graph. A {\em cyclet} of $\G$ is a sequence $(v_0, v_1, \ldots , v_{d-1})$ of pairwise distinct vertices of $\G$ such that $v_i v_{i+1}$ is an edge of $\G$ for each $i \in \ZZ_{d}$ (computation of indices is performed modulo $d$). The {\em length} of a cyclet is the length of the corresponding sequence. A cyclet can thus be thought of as a rooted oriented cycle and its length is the length of this {\em underlying cycle}. The {\em edges} of a cyclet are the edges of the corresponding underlying cycle. Therefore, a cyclet of length $d$ has precisely $d$ edges. We do not consider potential chordal edges as edges of the cyclet.

\begin{construction}
\label{con:general}
Let $\G$ and $\G'$ be graphs and let $d \geq 3$ be an integer. Let $C = (u_0, u_1, \ldots , u_{d-1})$ and $C' = (v_0, v_1, \ldots , v_{d - 1})$ be cyclets in $\G$ and $\G'$, respectively. Then the {\em merge} of $\G$ and $\G'$ {\em with respect} to $C$ and $C'$ is the graph $\G\oplus^{C}_{C'} \G'$, obtained by taking the disjoint union of $\G$ and $\G'$, deleting the $d$ edges of each of $C$ and $C'$ (we keep the vertices as well as potential chordal edges of the underlying cycles), and adding the edges $u_i v_{i+1}$ and $v_i u_{i+1}$ for each $i \in \ZZ_{d}$ (subscripts are computed modulo $d$).
\end{construction} 

\begin{figure}[htbp]
\begin{center}
\subfigure{
\begin{tikzpicture}[scale = 1.2]
\foreach \j in {0,1,2}
{
\node[vtx, fill=black, inner sep = 2pt, label=270:{\footnotesize $x_{\j}$}] (x\j) at (1.6*\j, 3.5) {};
\node[vtx, fill=black, inner sep = 2pt, label=90:{\footnotesize $y_{\j}$}] (y\j) at (1.6*\j, 4.5) {};
\node[vtx, fill=black, inner sep = 2pt, label=270:{\footnotesize $x'_{\j}$}] (xx\j) at (1.6*\j, 1) {};
\node[vtx, fill=black, inner sep = 2pt, label=90:{\footnotesize $y'_{\j}$}] (yy\j) at (1.6*\j, 2) {};
}
\foreach \j in {0,1}
{
\pgfmathtruncatemacro{\jj}{\j+1}
\draw[thick] (x\j) -- (x\jj);
\draw[thick] (x\j) -- (y\jj);
\draw[thick] (y\j) -- (y\jj);
\draw[thick] (y\j) -- (x\jj);
\draw[thick] (xx\j) -- (xx\jj);
\draw[thick] (xx\j) -- (yy\jj);
\draw[thick] (yy\j) -- (yy\jj);
\draw[thick] (yy\j) -- (xx\jj);
}
\draw[thick] (x0) -- (y2);
\draw[thick] (xx0) -- (yy2);
\draw[thick] (y0) -- (x2);
\draw[thick] (yy0) -- (xx2);
\draw[thick] (x0) to[bend left = 15] (x2);
\draw[thick] (xx0) to[bend left = 15] (xx2);
\draw[thick] (y0) to[bend right = 15] (y2);
\draw[thick] (yy0) to[bend right = 15] (yy2);

\draw[opacity = .5, blue, line width = 3 pt] (x0) -- (x1) -- (y2) -- (y1) -- (x0);
\draw[opacity = .5, blue, line width = 3 pt] (xx0) -- (xx1) -- (yy2) -- (yy1) -- (xx0);
\draw[opacity = .3, dashed, thick] (x0) -- (xx1) -- (y2) -- (yy1) -- (x0);
\draw[opacity = .3, dashed, thick] (xx0) -- (x1) -- (yy2) -- (y1) -- (xx0);
\end{tikzpicture}
}
\hspace{8mm}
\subfigure{
\begin{tikzpicture}[scale = 1.2]
\node[vtx, fill=black, inner sep = 2pt, label=300:{\footnotesize $x_0$}] (x0) at (120 : 1) {};
\node[vtx, fill=black, inner sep = 2pt, label=120:{\footnotesize $y_2$}] (y2) at (120 : 2) {};
\node[vtx, fill=black, inner sep = 2pt, label=240:{\footnotesize $x_2$}] (x2) at (60 : 1) {};
\node[vtx, fill=black, inner sep = 2pt, label=60:{\footnotesize $y_0$}] (y0) at (60 : 2) {};
\node[vtx, fill=black, inner sep = 2pt, label=180:{\footnotesize $x_1$}] (x1) at (0 : 1) {};
\node[vtx, fill=black, inner sep = 2pt, label=0:{\footnotesize $y_1$}] (y1) at (0 : 2) {};
\node[vtx, fill=black, inner sep = 2pt, label=120:{\footnotesize $x'_0$}] (xx0) at (300 : 1) {};
\node[vtx, fill=black, inner sep = 2pt, label=300:{\footnotesize $y'_2$}] (yy2) at (300 : 2) {};
\node[vtx, fill=black, inner sep = 2pt, label=60:{\footnotesize $x'_2$}] (xx2) at (240 : 1) {};
\node[vtx, fill=black, inner sep = 2pt, label=240:{\footnotesize $y'_0$}] (yy0) at (240 : 2) {};
\node[vtx, fill=black, inner sep = 2pt, label=0:{\footnotesize $x'_1$}] (xx1) at (180 : 1) {};
\node[vtx, fill=black, inner sep = 2pt, label=180:{\footnotesize $y'_1$}] (yy1) at (180 : 2) {};
\draw[thick] (x0) -- (x2) -- (y0) -- (y1) -- (x2) -- (x1) -- (y0) -- (y2) -- (x0);
\draw[thick] (xx0) -- (xx2) -- (yy0) -- (yy1) -- (xx2) -- (xx1) -- (yy0) -- (yy2) -- (xx0);
\draw[thick] (x0) -- (xx1) -- (y2) -- (yy1) -- (x0);
\draw[thick] (xx0) -- (x1) -- (yy2) -- (y1) -- (xx0);
\draw[opacity = .5, magenta, line width = 3 pt] (x0) -- (xx1) -- (y2) -- (yy1) -- (x0);
\draw[opacity = .5, magenta, line width = 3 pt] (xx0) -- (x1) -- (yy2) -- (y1) -- (xx0);
\end{tikzpicture}
}
\caption{The merged graph $\G\oplus^{C}_{C'} \G'$ where $\G = \G' = \Wr(3)$ and $C = C' = (x_0,x_1,y_2,y_1)$.}
\label{fig:merged_example}
\end{center}
\end{figure}
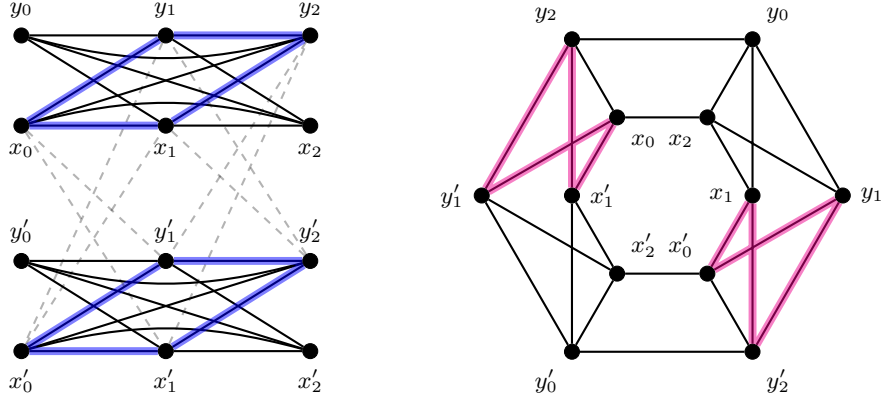

Observe that in the resulting merged graph $\G\oplus^{C}_{C'} \G'$, each vertex has the same valence as it had in its corresponding graph $\G$ or $\G'$, since we first delete two incident edges of a given vertex of $C$ (or $C'$) and then also add two edges incident to it. Note also that the added edges between the vertices of $C$ and $C'$ give rise to a cycle of length $2d$ when $d$ is odd and give rise to two vertex-disjoint cycles of length $d$ when $d$ is even. In all of our examples in this paper $d$ will in fact always be even. For example, if we let $\G = \G' = \Wr(3)$ and set $C = C' = (x_0,x_1,y_2,y_1)$, then the resulting merged graph $\G\oplus^{C}_{C'} \G'$ turns out to be the quasi-wreath graph $\QWr(3,3)$, as can be observed on Figure~\ref{fig:merged_example}, where on the left side of the figure the disjoint union of $\G$ and $\G'$ is presented, the cyclets $C$ and $C'$ are highlighted in blue, and the edges that will be added are indicated with dashed lines. On the right side of the figure the resulting merged graph $\G\oplus^{C}_{C'} \G'$ is presented and the new edges are highligted in magenta. The vertices are repositioned to show that the resulting graph is indeed $\QWr(3,3)$.

We now show that when $\G$ and $\G'$ are regular distance magic graphs and the corresponding distance magic labelings have certain nice enough properties with respect to the cyclets $C$ and $C'$, the merged graph $\G\oplus^{C}_{C'} \G'$ is also distance magic. Before proceeding we introduce the following additional terminology. 

Suppose $\G$ is a graph and let $\cP = \{A,B\}$ be any bipartition of its vertex set (where endvertices of an edge may belong to the same part). An edge $uv$ of $\G$ is said to be a {\em link with respect to $\cP$} whenever one of $u$ and $v$ is in $A$ and the other is in $B$, and is said to be a {\em non-link with respect to $\cP$} otherwise. A cycle (or cyclet) of even length of $\G$ is then said to be {\em alternating with respect to $\cP$} if its edges alternate between being links and non-links. For brevity we make an agreement that whenever the bipartition $\cP$ is clear from the context we simply speak of links, non-links and alternating cycles and cyclets. 

Now, suppose $\G$ is a graph such that each of its vertices has even valence. Then a bipartition $\cP = \{A, B\}$ of the vertex set of $\G$ is said to be {\em balanced} if each vertex of $\G$ has the same number of neighbors in $A$ as it does in $B$. In other words, for each vertex half of the edges incident to it are links while the other half are non-links. We mention that what we call a balanced bipartition in this paper is called a $2$-partition in~\cite{CicNik17} (defined only for regular graphs). 

Observe that if $\G$ is a regular distance magic graph and $\ell$ is a corresponding distance magic labeling, then one has a natural bipartition of the vertex set (not necessarily balanced) by letting $A$ consist of all the vertices $v$ with $\ell(v) \geq 0$ and $B$ of the remaining ones. Throughout the paper, when we speak of links and non-links with respect to $\ell$ and of alternating cycles and cyclets with respect to $\ell$, we are referring to this natural bipartition. As before, if $\ell$ is clear from the context we simply speak of links, non-links and alternating cycles and cyclets. Observe that the links are precisely the edges $uv$ for which precisely one of $\ell(u)$ and $\ell(v)$ is negative. 

\begin{proposition}
\label{pro:con}
Let $k \geq 2$ be an integer and let $\G$ and $\G'$ be $(2k)$-regular distance magic graphs with corresponding distance magic labelings $\ell$ and $\ell'$. Let $C$ and $C'$ be cyclets of even length $d$ as in Construction~\ref{con:general}. With notation from Construction~\ref{con:general} suppose each of the following holds:
\begin{itemize}
\item[(i)] The bipartition of $V(\G')$ corresponding to $\ell'$ is balanced.
\item[(ii)] The cyclet $C'$ is alternating with respect to $\ell'$.
\item[(iii)] $\ell(u_{i-1}) + \ell(u_{i+1}) = \ell'(v_{i-1}) + \ell'(v_{i+1})$ holds for each $i \in \ZZ_{d}$.
\end{itemize}
Then the merged graph $\G\oplus^{C}_{C'} \G'$ is a $(2k)$-regular distance magic graph.
\end{proposition}

\begin{proof}
Let $n$ and $n'$ be the orders of the graphs $\G$ and $\G'$, respectively. Let $\cP' = \{A',B'\}$ be the bipartition of $V(\G')$ corresponding to $\ell'$, where $A' = \{v' \in V(\G') \colon \ell'(v') \geq 0\}$. By assumption (i), this bipartition is balanced, and thus a simple double counting argument shows that $|A'| = |B'|$, implying that $n'$ is even. Set $\tilde{\G} = \G\oplus^{C}_{C'} \G'$ and let $\tilde{\ell}$ be the labeling of its vertices defined as follows:
 \begin{equation}
 \label{eq:ell_tilde}
 	\tilde{\ell}(w) = \left\{\begin{array}{ccc}
 		\ell(w) & : & w \in V(\G),\\
 		\ell'(w) + n & : & w \in A', \\
 		\ell'(w) - n & : & w \in B'.\end{array}\right.
 \end{equation}
We claim that $\tilde{\ell}$ is a distance magic labeling of $\tilde{\G}$. That $\tilde{\ell}$ is a bijection from $V(\tilde{\G})$ to the set $\mathcal{I}_{n+n'}$, follows from the definition of the sets $A'$ and $B'$, and the fact that $\ell$ and $\ell'$ are distance magic labelings of $\G$ and $\G'$, respectively. 
 
Consider first a vertex $u \in V(\G)$ and recall that $\sum_{w \in \G(u)}\ell(w) = 0$ since $\ell$ is distance magic. If $u$ is not a vertex of the cyclet $C$, then $\tilde{\G}(u) = \G(u)$, and so \eqref{eq:ell_tilde} implies that $\sum_{w \in \tilde{\G}(u)}\tilde{\ell}(w) = \sum_{w \in \G(u)}\ell(w) = 0$. If however $u = u_i$ for some $i \in \ZZ_d$, then $\tilde{\G}(u) = (\G(u) \setminus \{u_{i-1}, u_{i+1}\}) \cup \{v_{i-1}, v_{i+1}\}$. By assumption (ii), the cyclet $C'$ is alternating with respect to $\ell'$, and so precisely one of $v_{i-1}$ and $v_{i+1}$ is in $A'$ (while the other is in $B'$). Due to assumption (iii), the definition of $\tilde{\ell}$ from~\eqref{eq:ell_tilde} yields
$$
	\sum_{w \in \tilde{\G}(u)}\tilde{\ell}(w) = \left(\sum_{w \in \G(u)}\ell(w)\right) - \ell(u_{i-1}) - \ell(u_{i+1}) + \ell'(v_{i-1}) + \ell'(v_{i+1}) = \sum_{w \in \G(u)}\ell(w) = 0.
$$

Consider now a vertex $v \in V(\G')$ and recall that $\sum_{w \in \G'(v)}\ell'(w) = 0$. By assumption (i), half of the neighbors of $v$ in $\G'$ belong to $A'$ and the other half to $B'$. Therefore, if $v$ is not a vertex of $C'$ (in which case $\tilde{\G}(v) = \G'(v)$), then~\eqref{eq:ell_tilde} implies $\sum_{w \in \tilde{\G}(v)}\tilde{\ell}(w) = \sum_{w \in \G(v)}\ell'(w) = 0$. Suppose then that $v = v_i$ for some $i \in \ZZ_d$. By assumption (ii), precisely one of $v_{i-1}$ and $v_{i+1}$ is in $A'$ (while the other is in $B'$),
and so assumption (i) implies that $|\tilde{\G}(v) \cap A'| = k-1 = |\tilde{\G}(v) \cap B'|$. Assumption (iii), together with~\eqref{eq:ell_tilde}, thus yields
$$
	\sum_{w \in \tilde{\G}(v)}\tilde{\ell}(w) = \ell(u_{i-1}) + \ell(u_{i+1}) + \sum_{w \in \tilde{\G}(v) \cap (A' \cup B')}\ell'(w) = \sum_{w \in \G'(v)}\ell'(w) = 0.
$$
This thus shows that $\tilde{\ell}$ is indeed a distance magic labeling of $\tilde{\G}$.
\end{proof}

The next corollary shows that if the distance magic labelings of $\G$ and $\G'$ are self-reverse and the cyclets $C$ and $C'$ are of particular kind, then the resulting distance magic labeling $\tilde{\ell}$ from the above proof is also self-reverse. 

\begin{corollary}
\label{cor:SRcon}
Let $k \geq 2$ and let $\G$ and $\G'$ be $(2k)$-regular distance magic graphs admitting self-reverse distance magic labelings $\ell$ and $\ell'$, respectively. Let $C = (u_0, u_1, \ldots , u_{d-1})$ and $C' = (v_0, v_1, \ldots , v_{d - 1})$ be cyclets of even length $d = 2d_0$ as in Construction~\ref{con:general}, and suppose the conditions (i) -- (iii) of Proposition~\ref{pro:con} all hold. Suppose that one of the following holds:
\begin{itemize}
\item[(i)] $u_{i+d_0} = u_i^\ell$ and $v_{i+d_0} = v_i^{\ell'}$ for all $i \in \ZZ_d$ (computations performed modulo $d$), or
\item[(ii)] $u_{i+d_0} = u_{d_0-1-i}^\ell$ and $v_{i+d_0} = v_{d_0-1-i}^{\ell'}$ for all $i \in \ZZ_d$ (computations performed modulo $d$).
\end{itemize}  
Then the labeling $\tilde{\ell}$ from the proof of Proposition~\ref{pro:con} is a distance magic self-reverse labeling of $\G\oplus^{C}_{C'} \G'$.
\end{corollary}

\begin{proof}
By Proposition~\ref{pro:con}, the labeling $\tilde{\ell}$ from its proof is a distance magic labeling of $\tilde{\G}$, where $\tilde{\G} = \G\oplus^{C}_{C'} \G'$. We thus only need to prove that it is self-reverse. Before proceeding observe that each of the assumptions (i) and (ii) of this corollary implies that if $u$ and $w$ are adjacent vertices of $\G$, then $u w$ is an edge of $C$ if and only if $u^\ell w^\ell$ is an edge of $C$. Similarly, if $v$ and $w$ are adjacent vertices of $\G'$, then $v w$ is an edge of $C'$ if and only if $v^{\ell'} w^{\ell'}$ is an edge of $C'$. Observe also that by the definition of $\tilde{\ell}$ and of the bipartition $\mathcal{P}' = \{A', B'\}$ from the proof of Proposition~\ref{pro:con}, it follows that $u^{\tilde{\ell}} = u^\ell$ for each $u \in V(\G)$ and $v^{\tilde{\ell}} = v^{\ell'}$ for each $v \in V(\G')$.

Suppose first that some $u, w \in V(\G)$ are adjacent in $\tilde{\G}$. By construction of $\tilde{\G}$, the edge $u w$ is not an edge of $C$, and so $u^\ell w^\ell$ is also not an edge of $C$. But since $\ell$ is self-reverse, $u^\ell w^\ell$ is an edge of $\G$ and thus also of $\tilde{\G}$. Since $u^{\tilde{\ell}} = u^\ell$ and $w^{\tilde{\ell}} = w^\ell$, this thus shows that $u^{\tilde{\ell}}$ is adjacent to $w^{\tilde{\ell}}$ in $\tilde{\G}$. 

The completely analogous argument that if $v,w \in V(\G')$ are adjacent in $\tilde{\G}$, then $v^{\tilde{\ell}}$ and $w^{\tilde{\ell}}$ are adjacent as well, is left to the reader. 

Suppose finally that $u \in V(\G)$ and $v \in V(\G')$ are adjacent in $\tilde{\G}$. Then $u = u_i$ and $v \in \{v_{i-1}, v_{i+1}\}$ for some $i \in \ZZ_d$. If the assumption (i) of this corollary holds, then 
$$
	u^{\tilde{\ell}} = u^\ell = u_{i+d_0} \quad \text{and} \quad v^{\tilde{\ell}} = v^{\ell'} \in \{v_{i-1}^{\ell'}, v_{i+1}^{\ell'}\} = \{v_{i+d_0 - 1}, v_{i+d_0 + 1}\},
$$
showing that $u^{\tilde{\ell}}$ and $v^{\tilde{\ell}}$ are adjacent in $\tilde{\G}$. If however, the assumption (ii) of this corollary holds, then
$$
	u^{\tilde{\ell}} = u^\ell = u_{2d_0-1-i} \quad \text{and} \quad v^{\tilde{\ell}} = v^{\ell'} \in \{v_{i-1}^{\ell'}, v_{i+1}^{\ell'}\} = \{v_{2d_0 - i}, v_{2d_0 -2 - i}\},
$$
again showing that $u^{\tilde{\ell}}$ and $v^{\tilde{\ell}}$ are adjacent in $\tilde{\G}$. The distance magic labeling $\tilde{\ell}$ is thus indeed self-reverse.
\end{proof}
 
A few remarks are in order. It can be verified that in order to fulfill the conditions of Proposition~\ref{pro:con} and Corollary~\ref{cor:SRcon}, the length $d$ of the cyclets $C$ and $C'$ must be divisible by $4$, where in the case that the condition (i) from Corollary~\ref{cor:SRcon} holds, we in addition require $d \equiv 4 \pmod{8}$ and $d \geq 12$. 

For our purposes the construction from Construction~\ref{con:general}, together with Corollary~\ref{cor:SRcon}, will suffice, but we point out the following possibility for a generalization. Namely, Construction~\ref{con:general} could be generalized in the sense that one would take a collection $\mathcal{C} = \{C_1, C_2, \ldots , C_m\}$ of cyclets from $\G$ and a collection $\mathcal{C}' = \{C'_1, C'_2, \ldots , C'_m\}$ of cyclets from $\G'$ such that each $C_i$ is of the same length as $C'_i$, and then construct a corresponding merged graph, where for each $i$ we delete the edges of $C_i$ and $C'_i$ and insert the new edges connecting the vertices of these two cyclets as in Construction~\ref{con:general}. If for instance no two cyclets in any of the two collections share a vertex, then the vertices in the resulting merged graph will again for sure retain their valence. Moreover, one can then generalize Proposition~\ref{pro:con}, where we require that all the cyclets from $\mathcal{C}'$ are alternating with respect to $\ell'$ and that each $C_i$ is compatible with $C'_i$ in the sense of condition (iii) from that proposition. This can then for instance be used to generalize Corollary~\ref{cor:SRcon} in the sense that we can take a cyclet $C_1$ of $\G$ and $C_1'$ of $\G'$ (of the same length as $C_1$), each of which visits at most one vertex from each pair $\{u, u^\ell\}$ and $\{v, v^{\ell'}\}$, and then let $C_2$ and $C'_2$ be the cyclets obtained from $C_1$ and $C_1'$, respectively, by exchanging each vertex $u$ of $C_1$ by $u^\ell$ and each vertex $v$ of $C'_1$ by $v^{\ell'}$. The above mentioned construction (where in the case that $\G$ is of odd order we require that $C_1$ avoids the central vertex) with $\mathcal{C} = \{C_1, C_2\}$ and $\mathcal{C'} = \{C'_1, C'_2\}$ then again yields a graph admitting a self-reverse distance magic labeling (provided of course that all the assumptions from Corollary~\ref{cor:SRcon}, except for the conditions (i) and (ii), hold).

Before illustrating Corollary~\ref{cor:SRcon} with a concrete example we also point out the following. When searching for suitable distance magic graphs with appropriate non-degenerate self-reverse distance magic labelings in which to search for cyclets that will enable us to apply Corollary~\ref{cor:SRcon}, it suffices to work with the corresponding quotients. What we need is either an alternating simple cycle (in the sense that the solid and dashed edges alternate on it) with an odd number of dashed edges (to fulfill condition (i) from Corollary~\ref{cor:SRcon}), or an alternating path of odd length, starting and ending with a solid edge, and such that the initial and terminal vertex each has a semiedge. Of course, in both cases we need to avoid the central vertex in the case of graphs of odd order to actually get a cycle. For instance, if we let $\G$ and $\ell$ be as in Figure~\ref{fig:Odd21quo}, then a possible suitable cyclet $C$ could be the cyclet of length $12$, corresponding to the simple cycle $(8,12,20,14,10,18)$ in the quotient (where the numbers correspond to the labels of the vertices of the quotient on the right part of Figure~\ref{fig:Odd21quo}). 

\begin{figure}[htbp]
\begin{center}
\subfigure{
\begin{tikzpicture}[scale = .4]
\node[vtx, fill=black, inner sep = 2pt, label=300:{$3$}] (A0) at (60: 4) {};
\node[vtx, fill=black, inner sep = 2pt, label=60:{$7$}] (A1) at (300: 4) {};
\node[vtx, fill=black, inner sep = 2pt, label=90:{$5$}] (A2) at (180: 4) {};
\node[vtx, fill=black, inner sep = 2pt, label=120:{$1$}] (A3) at (0:0) {};
\draw[thick] (A0) -- (A1);
\draw[thick, opacity = .5, blue, line width = 3 pt] (A0) -- (A1);
\draw[dashed, thick] (A0) -- (A2);
\draw[thick] (A0) -- (A3);
\draw[thick] (A1) -- (A2);
\draw[dashed, thick] (A1) -- (A3);
\draw[thick] (A2) -- (A3);

\begin{scope}[on background layer]
\draw[dashed, thick, shorten <= -5mm] (A0.60) -- (A0.240);
\draw[dashed, thick, shorten <= -5mm] (A1.300) -- (A1.120);
\draw[thick, opacity = .5, blue, line width = 3 pt, shorten <= -5mm] (A0.60) -- (A0.240);
\draw[thick, opacity = .5, blue, line width = 3 pt, shorten <= -5mm] (A1.300) -- (A1.120);
\draw[dashed, thick, shorten <= -5mm] (A2.180) -- (A2.0);
\draw[dashed, thick, shorten <= -5mm] (A3.0) -- (A3.180);
\end{scope}
\end{tikzpicture}
}
\hspace{-3mm}
\subfigure{
\begin{tikzpicture}[scale = .4]
\node[vtx, fill=black, inner sep = 2pt, label=270:{$3$}] (A0) at (0: 4) {};
\node[vtx, fill=black, inner sep = 2pt, label=240:{$1$}] (A1) at (120: 4) {};
\node[vtx, fill=black, inner sep = 2pt, label=120:{$5$}] (A2) at (240: 4) {};
\node[vtx, fill=black, inner sep = 2pt, label=300:{$7$}] (A3) at (0:0) {};
\draw[thick] (A0) -- (A1);
\draw[dashed, thick] (A0) -- (A2);
\draw[thick] (A0) -- (A3);
\draw[thick] (A1) -- (A2);
\draw[thick, opacity = .5, blue, line width = 3 pt] (A1) -- (A2);
\draw[dashed, thick] (A1) -- (A3);
\draw[thick] (A2) -- (A3);

\begin{scope}[on background layer]
\draw[dashed, thick, shorten <= -5mm] (A0.0) -- (A0.180);
\draw[dashed, thick, shorten <= -5mm] (A1.120) -- (A1.300);
\draw[thick, opacity = .5, blue, line width = 3 pt, shorten <= -5mm] (A1.120) -- (A1.300);
\draw[dashed, thick, shorten <= -5mm] (A2.240) -- (A2.60);
\draw[thick, opacity = .5, blue, line width = 3 pt, shorten <= -5mm] (A2.240) -- (A2.60);
\draw[dashed, thick, shorten <= -5mm] (A3.180) -- (A3.0);
\end{scope}
\end{tikzpicture}
}
\hspace{8mm}
\subfigure{
\begin{tikzpicture}[scale = .4]
\node[vtx, fill=black, inner sep = 2pt, label=270:{$3$}] (A0) at (60: 4) {};
\node[vtx, fill=black, inner sep = 2pt, label=90:{$7$}] (A1) at (300: 4) {};
\node[vtx, fill=black, inner sep = 2pt, label=90:{$5$}] (A2) at (180: 4) {};
\node[vtx, fill=black, inner sep = 2pt, label=120:{$1$}] (A3) at (0:0) {};
\draw[dashed, thick] (A0) -- (A2);
\draw[thick] (A0) -- (A3);
\draw[thick] (A1) -- (A2);
\draw[dashed, thick] (A1) -- (A3);
\draw[thick] (A2) -- (A3);

\begin{scope}[on background layer]
\draw[dashed, thick, shorten <= -5mm] (A2.180) -- (A2.0);
\draw[dashed, thick, shorten <= -5mm] (A3.0) -- (A3.180);
\end{scope}

\phantom{
\foreach \j in {0,1,2}{
\pgfmathtruncatemacro{\k}{2*\j+3}
\node[vtx, fill=black, inner sep = 2pt, label=270-\j*60:{$\k$}] (B\j) at (360*\j/3: 4) {};
}
\node[vtx, fill=black, inner sep = 2pt, label=300:{$1$}] (B3) at (0:0) {};
}
\node[vtx, fill=black, inner sep = 2pt, label=90:{$11$}] (C0) at ($(B0) + (7,0)$){};
\node[vtx, fill=black, inner sep = 2pt, label=270:{$9$}] (C1) at ($(B1) + (7,0)$){};
\node[vtx, fill=black, inner sep = 2pt, label=90:{$13$}] (C2) at ($(B2) + (7,0)$){};
\node[vtx, fill=black, inner sep = 2pt, label=60:{$15$}] (C3) at ($(B3) + (7,0)$){};

\draw[thick] (C0) -- (C3);
\draw[dashed, thick] (C1) -- (C3);
\draw[dashed, thick] (C0) -- (C2);
\draw[thick] (C0) -- (C1);
\draw[thick] (C2) -- (C3);
\draw[dashed, thick] (A0) -- (C1);
\draw[dashed, thick] (A1) -- (C2);
\draw[thick] (A0) -- (C2);
\draw[thick] (A1) -- (C1);
\begin{scope}[on background layer]
\draw[dashed, thick, shorten <= -5mm] (C0.0) -- (C0.180);
\draw[dashed, thick, shorten <= -5mm] (C3.180) -- (C3.0);
\draw[thick, opacity = .5, magenta, line width = 3 pt] (A0) -- (C1);
\draw[thick, opacity = .5, magenta, line width = 3 pt] (A0) -- (C2);
\draw[thick, opacity = .5, magenta, line width = 3 pt] (A1) -- (C1);
\draw[thick, opacity = .5, magenta, line width = 3 pt] (A1) -- (C2);
\end{scope}

\end{tikzpicture}
}
\caption{A (non-degenerate) self-reverse distance magic labeling of $\Wr(8)$ using Corollary~\ref{cor:SRcon}.}
\label{fig:SR16}
\end{center}
\end{figure}
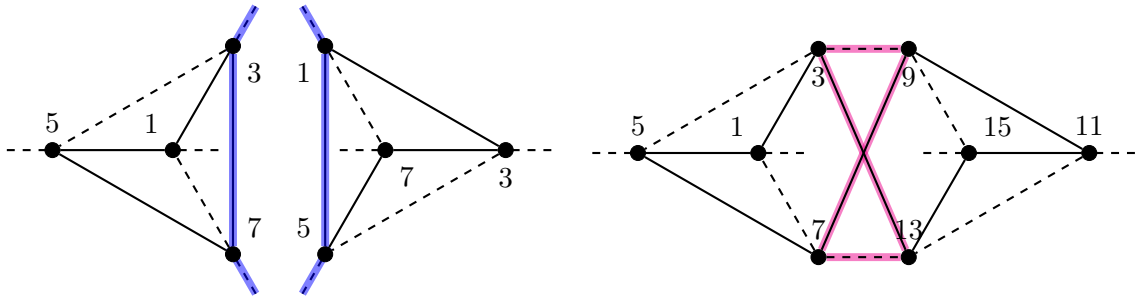

In all of our examples, however, we will be using the most simple variant of the above situation, where we simply take a $4$-cycle corresponding to a solid edge whose endvertices both have semiedges in the corresponding quotient. Observe that in this case the condition (iii) from Proposition~\ref{pro:con} will be satisfied if and only if the  difference of the two labels in the corresponding quotient will be the same for both choices (in $\G$ and $\G'$). For instance, if we let $\G = \G'$ be the graph $\Wr(4)$, let $\ell = \ell'$ be the labeling from Figure~\ref{fig:W(4)quo}, and let $C$ be the cyclet of length $4$, corresponding to the ordered pair $(3,7)$ in the quotient (where we interpret this as the cyclet visiting the vertices with labels $3$, $7$, $-7$ and $-3$ in this order), and $C'$ the cyclet of length $4$, corresponding to the ordered pair $(1,5)$ in the quotient, then the resulting merged graph $\G\oplus^{C}_{C'} \G'$ is the wreath graph $\Wr(8)$ of order $16$ and the resulting (non-degenerate) self-reverse distance magic labeling is the one presented (via the corresponding quotient) in Figure~\ref{fig:SR16}. 

\section{The classification}
\label{sec:classification}

We are now ready to determine the set of all positive integers $n$ such that a connected tetravalent graph of order $n$ admitting a self-reverse distance magic labeling exists. Since we already know that for each $m \geq 3$ the wreath graph $\Wr(m)$ (which is of order $2m$) admits a (degenerate) self-reverse distance magic labeling, it seems that the even integers are not interesting. However, we can of course also ask about the existence of non-degenerate distance magic labelings and we can even require that the underlying graph is not a wreath graph.

Using a personal computer one can determine all non-degenerate self-reverse distance magic labelings for connected tetravalent graphs up to some small enough order. Our computations reveal (but see also~\cite{RozSpa25}, where all the tetravalent distance magic graphs up to order $16$ were given) that up to order $16$ the only connected tetravalent graphs admitting a self-reverse distance magic labeling are the wreath graphs $\Wr(m)$, where $3 \leq m \leq 8$. Moreover, by Proposition~\ref{pro:wreath}, the only examples among these admitting a non-degenerate self-reverse distance magic labeling are $\Wr(4)$ and $\Wr(8)$ (see Figures~\ref{fig:W(4)quo} and~\ref{fig:SR16} for concrete examples of such labelings).

We were in fact able to compute all non-degenerate self-reverse distance magic labelings for connected tetravalent graphs up to order $30$, up to equivalence (in the sense of the definition from the last paragraph of Section~\ref{sec:Self}). In Table~\ref{tab:data} we provide the obtained data, where for each positive integer $n$ with $16 \leq n \leq 30$ we state the number of nonequivalent non-degenerate self-reverse distance magic labelings ($\# \text{SR}$) of connected tetravalent graphs of order $n$, together with the number of pairwise nonisomorphic graphs ($\# \text{gr}$) that these labelings give rise to. We explain the meaning of the last row of the table in the next section.
\begin{table}
$$
\begin{array}{|c||@{\ }c@{\ }|@{\ }c@{\ }|@{\ }c@{\ }|@{\ }c@{\ }|@{\ }c@{\ }|@{\ }c@{\ }|@{\ }c@{\ }|@{\ }c@{\ }|@{\ }c@{\ }|@{\ }c@{\ }|@{\ }c@{\ }|@{\ }c@{\ }|@{\ }c@{\ }|@{\ }c@{\ }|@{\ }c@{\ }|}
\hline
n & 16 & 17 & 18 & 19 & 20 & 21 & 22 & 23 & 24 & 25 & 26 & 27 & 28 & 29 & 30 \\ \hline \hline
\# \text{SR} & 48 & 0 & 136 & 0 & 66 & 57 & 0 & 675 & 11156 & 3063 & 31562 & 10951 & 35402 & 68837 & 229716 \\
\hline 
\# \text{gr} & 1 & 0 & 2 & 0 & 2 & 7 & 0 & 80 & 9 & 522 & 37 & 2647 & 342 & 22893 & 4151 \\
\hline
\# \text{VT} & 1 & 0 & 1 & 0 & 1 & 0 & 0 & 0 & 3 & 0 & 0 & 0 & 0 & 0 & 1 \\
\hline
\end{array}
$$
\caption{Nonequivalent non-degenerate self-reverse distance magic labelings of connected tetravalent graphs, the number of corresponding graphs and the number of vertex-transitive ones.}
\label{tab:data}
\end{table}

The data suggests that for each $n \geq 23$, at least one connected tetravalent graph admitting a non-degenerate self-reverse distance magic labeling exists. We show that this is indeed the case.

\begin{theorem}
\label{the:classification}
Let $n \geq 5$ be an integer. Then the following all hold:
\begin{itemize}
\itemsep = 0pt
\item[(i)] There exists a connected tetravalent graph of order $n$ admitting a self-reverse distance magic labeling if and only if either $n$ is even with $n \geq 6$, or $n$ is odd with $n \geq 21$.
\item[(ii)] There exists a connected tetravalent graph of order $n$ admitting a non-degenerate self-reverse distance magic labeling if and only if either $n \geq 23$ or $n \in \{8, 16, 18, 20, 21\}$. 
\item[(iii)] There exists a connected tetravalent graph of order $n$, which is not a wreath graph, but admits a self-reverse distance magic labeling, if and only if $n \geq 18$ with $n \notin \{19, 22\}$.
\end{itemize}
\end{theorem}

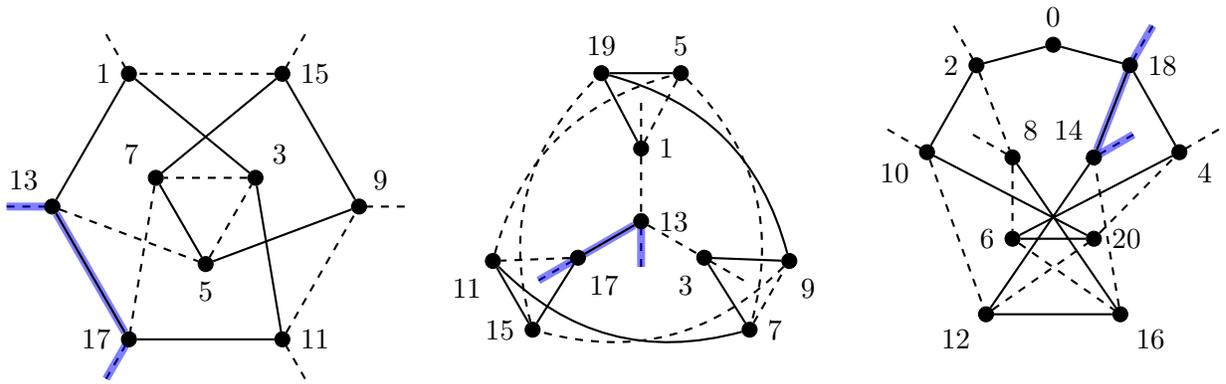
\begin{figure}[htbp]
\begin{center}
\subfigure{
\begin{tikzpicture}[scale = 0.51]
\node[vtx, fill=black, inner sep = 2pt, label=60:{$9$}] (A0) at (0: 4) {};
\node[vtx, fill=black, inner sep = 2pt, label=0:{$15$}] (A1) at (60: 4) {};
\node[vtx, fill=black, inner sep = 2pt, label=180:{$1$}] (A2) at (120: 4) {};
\node[vtx, fill=black, inner sep = 2pt, label=120:{$13$}] (A3) at (180: 4) {};
\node[vtx, fill=black, inner sep = 2pt, label=180:{$17$}] (A4) at (240: 4) {};
\node[vtx, fill=black, inner sep = 2pt, label=0:{$11$}] (A5) at (300: 4) {};
\node[vtx, fill=black, inner sep = 2pt, label=30:{$3$}] (B0) at (30: 1.5) {};
\node[vtx, fill=black, inner sep = 2pt, label=150:{$7$}] (B1) at (150: 1.5) {};
\node[vtx, fill=black, inner sep = 2pt, label=270:{$5$}] (B2) at (270: 1.5) {};
\draw[thick] (A0) -- (A1) -- (B1) -- (B2) -- (A0);
\draw[thick] (A2) -- (A3) -- (A4) -- (A5) -- (B0) -- (A2);
\draw[dashed, thick] (A1) -- (A2);
\draw[dashed, thick] (A5) -- (A0);
\draw[dashed, thick] (A3) -- (B2) -- (B0) -- (B1) -- (A4);
\begin{scope}[on background layer]
\draw[dashed, thick, shorten <= -5mm] (A0.0) -- (A0.180);
\draw[dashed, thick, shorten <= -5mm] (A1.60) -- (A1.240);
\draw[dashed, thick, shorten <= -5mm] (A2.120) -- (A2.300);
\draw[dashed, thick, shorten <= -5mm] (A3.180) -- (A3.0);
\draw[dashed, thick, shorten <= -5mm] (A4.240) -- (A4.60);
\draw[dashed, thick, shorten <= -5mm] (A5.300) -- (A5.120);
\draw[thick, opacity = .5, blue, line width = 3 pt] (A3) -- (A4);
\draw[thick, shorten <= -5mm, opacity = .5, blue, line width = 3 pt] (A3.180) -- (A3.0);
\draw[thick, shorten <= -5mm, opacity = .5, blue, line width = 3 pt] (A4.240) -- (A4.60);
\end{scope}
\end{tikzpicture}
}
\hspace{1mm}
\subfigure{
\begin{tikzpicture}[scale = 0.51]
\node[vtx, fill=black, inner sep = 2pt, label=0:{$13$}] (A0) at (0: 0) {};
\node[vtx, fill=black, inner sep = 2pt, label=0:{$1$}] (B0) at (90: 1.9) {};
\node[vtx, fill=black, inner sep = 2pt, label=280:{$17$}] (B1) at (210: 1.9) {};
\node[vtx, fill=black, inner sep = 2pt, label=260:{$3$}] (B2) at (330: 1.9) {};
\node[vtx, fill=black, inner sep = 2pt, label=90:{$5$}] (C0) at (75: 4) {};
\node[vtx, fill=black, inner sep = 2pt, label=90:{$19$}] (D0) at (105: 4) {};
\node[vtx, fill=black, inner sep = 2pt, label=260:{$11$}] (C1) at (195: 4) {};
\node[vtx, fill=black, inner sep = 2pt, label=180:{$15$}] (D1) at (225: 4) {};
\node[vtx, fill=black, inner sep = 2pt, label=0:{$7$}] (C2) at (315: 4) {};
\node[vtx, fill=black, inner sep = 2pt, label=280:{$9$}] (D2) at (345: 4) {};
%
\draw[thick] (B0) -- (D0) -- (C0);
\draw[thick] (A0) -- (B1) -- (D1) -- (C1);
\draw[thick] (C2) -- (B2) -- (D2);
\draw[dashed, thick] (C0) -- (B0) -- (A0) -- (B2);
\draw[dashed, thick] (B1) -- (C1);
\draw[dashed, thick] (C2) -- (D2);
\draw[thick] (D0) to [bend left = 30] (D2);
\draw[thick] (C2) to [bend left = 30] (C1);
\draw[dashed, thick] (C1) to [bend left = 30] (C0);
\draw[dashed, thick] (C0) to [bend left = 30] (C2);
\draw[dashed, thick] (D1) to [bend left = 30] (D0);
\draw[dashed, thick] (D2) to [bend left = 30] (D1);
\begin{scope}[on background layer]
\draw[dashed, thick, shorten <= -5mm] (A0.270) -- (A0.90);
\draw[dashed, thick, shorten <= -5mm] (B0.90) -- (B0.270);
\draw[dashed, thick, shorten <= -5mm] (B1.210) -- (B1.30);
\draw[dashed, thick, shorten <= -5mm] (B2.330) -- (B2.150);
\draw[thick, opacity = .5, blue, line width = 3 pt] (A0) -- (B1);
\draw[thick, shorten <= -5mm, opacity = .5, blue, line width = 3 pt] (A0.270) -- (A0.90);
\draw[thick, shorten <= -5mm, opacity = .5, blue, line width = 3 pt] (B1.210) -- (B1.30);
\end{scope}
\end{tikzpicture}
}
\hspace{1mm}
\subfigure{
\begin{tikzpicture}[scale = 0.51]
\node[vtx, fill=black, inner sep = 2pt, label=90:{$0$}] (C) at (90: 4) {};
\node[vtx, fill=black, inner sep = 2pt, label=95:{$14$}] (A0) at (45: 1.5) {};
\node[vtx, fill=black, inner sep = 2pt, label=85:{$8$}] (A1) at (135: 1.5) {};
\node[vtx, fill=black, inner sep = 2pt, label=180:{$6$}] (A2) at (225: 1.5) {};
\node[vtx, fill=black, inner sep = 2pt, label=0:{$20$}] (A3) at (315: 1.5) {};
\node[vtx, fill=black, inner sep = 2pt, label=350:{$4$}] (B0) at (20: 3.5) {};
\node[vtx, fill=black, inner sep = 2pt, label=190:{$10$}] (B1) at (160: 3.5) {};
\node[vtx, fill=black, inner sep = 2pt, label=225:{$12$}] (B2) at (240: 3.5) {};
\node[vtx, fill=black, inner sep = 2pt, label=315:{$16$}] (B3) at (300: 3.5) {};
\node[vtx, fill=black, inner sep = 2pt, label=0:{$18$}] (D0) at (60: 4) {};
\node[vtx, fill=black, inner sep = 2pt, label=180:{$2$}] (D1) at (120: 4) {};
\draw[thick] (A1) -- (B3) -- (B2) -- (A0) -- (D0) -- (C) -- (D1) -- (B1) -- (A3) -- (A2) -- (B0) -- (D0);
\draw[dashed, thick] (D1) -- (A1) -- (A2) -- (B3) -- (A0);
\draw[dashed, thick] (B1) -- (B2) -- (A3) -- (B0);
\begin{scope}[on background layer]
\draw[dashed, thick, shorten <= -5mm] (D1.120) -- (D1.300);
\draw[dashed, thick, shorten <= -5mm] (D0.60) -- (D0.240);
\draw[dashed, thick, shorten <= -5mm] (B1.150) -- (B1.330);
\draw[dashed, thick, shorten <= -5mm] (B0.30) -- (B0.210);
\draw[dashed, thick, shorten <= -5mm] (A1.150) -- (A1.330);
\draw[dashed, thick, shorten <= -5mm] (A0.30) -- (A0.210);
\draw[thick, opacity = .5, blue, line width = 3pt] (A0) -- (D0);
\draw[thick, shorten <= -5mm, opacity = .5, blue, line width = 3 pt] (D0.60) -- (D0.240);
\draw[thick, shorten <= -5mm, opacity = .5, blue, line width = 3 pt] (A0.30) -- (A0.210);
\end{scope}
\end{tikzpicture}
}
\caption{Examples of orders $18$, $20$ and $21$ admitting ``infinite'' extension.}
\label{fig:n18_20_21}
\end{center}
\end{figure}

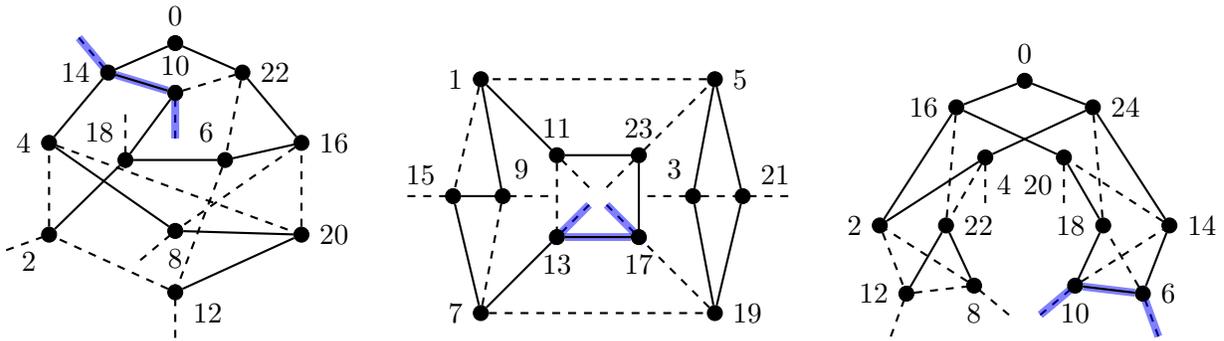
\begin{figure}[htbp]
\begin{center}
\subfigure{
\begin{tikzpicture}[scale = 0.51]
\node[vtx, fill=black, inner sep = 2pt, label=90:{$0$}] (C) at (90: 3.8) {};
\node[vtx, fill=black, inner sep = 2pt, label=100:{$6$}] (A0) at (30: 1.5) {};
\node[vtx, fill=black, inner sep = 2pt, label=100:{$18$}] (A1) at (150: 1.5) {};
\node[vtx, fill=black, inner sep = 2pt, label=270:{$8$}] (A2) at (270: 1.1) {};
\node[vtx, fill=black, inner sep = 2pt, label=0:{$16$}] (B0) at (20: 3.5) {};
\node[vtx, fill=black, inner sep = 2pt, label=0:{$22$}] (B1) at (60: 3.5) {};
\node[vtx, fill=black, inner sep = 2pt, label=180:{$14$}] (B2) at (120: 3.5) {};
\node[vtx, fill=black, inner sep = 2pt, label=180:{$4$}] (B3) at (160: 3.5) {};
\node[vtx, fill=black, inner sep = 2pt, label=250:{$2$}] (B4) at (200: 3.5) {};
\node[vtx, fill=black, inner sep = 2pt, label=350:{$12$}] (B5) at (270: 2.7) {};
\node[vtx, fill=black, inner sep = 2pt, label=0:{$20$}] (B6) at (340: 3.5) {};
\node[vtx, fill=black, inner sep = 2pt, label=90:{$10$}] (D) at (90: 2.5) {};
\draw[thick] (B4) -- (A1) -- (D) -- (B2) -- (B3) -- (A2) -- (B6) -- (B5);
\draw[thick] (A1) -- (A0) -- (B0) -- (B1) -- (C) -- (B2);
\draw[dashed, thick] (D) -- (B1) -- (A0) -- (B5) -- (B4) -- (B3) -- (B6) -- (B0) -- (A2);
\begin{scope}[on background layer]
\draw[dashed, thick, shorten <= -5mm] (A2.220) -- (A2.40);
\draw[dashed, thick, shorten <= -5mm] (A1.90) -- (A1.270);
\draw[dashed, thick, shorten <= -5mm] (D.270) -- (D.90);
\draw[dashed, thick, shorten <= -5mm] (B5.270) -- (B5.90);
\draw[dashed, thick, shorten <= -5mm] (B4.200) -- (B4.20);
\draw[dashed, thick, shorten <= -5mm] (B2.130) -- (B2.310);
\draw[thick, opacity = .5, blue, line width = 3pt] (B2) -- (D);
\draw[thick, shorten <= -5mm, opacity = .5, blue, line width = 3 pt] (D.270) -- (D.90);
\draw[thick, shorten <= -5mm, opacity = .5, blue, line width = 3 pt] (B2.130) -- (B2.310);
\end{scope}
\end{tikzpicture}
}
\hspace{0mm}
\subfigure{
\begin{tikzpicture}[scale = 0.55]
\node[vtx, fill=black, inner sep = 2pt, label=100:{$3$}] (A0) at (0: 2.3) {};
\node[vtx, fill=black, inner sep = 2pt, label=90:{$23$}] (A1) at (45: 1.4) {};
\node[vtx, fill=black, inner sep = 2pt, label=90:{$11$}] (A2) at (135: 1.4) {};
\node[vtx, fill=black, inner sep = 2pt, label=80:{$9$}] (A3) at (180: 2.3) {};
\node[vtx, fill=black, inner sep = 2pt, label=270:{$13$}] (A4) at (225: 1.4) {};
\node[vtx, fill=black, inner sep = 2pt, label=270:{$17$}] (A5) at (315: 1.4) {};
\node[vtx, fill=black, inner sep = 2pt, label=10:{$21$}] (B0) at (0: 3.5) {};
\node[vtx, fill=black, inner sep = 2pt, label=0:{$5$}] (B1) at (45: 4) {};
\node[vtx, fill=black, inner sep = 2pt, label=180:{$1$}] (B2) at (135: 4) {};
\node[vtx, fill=black, inner sep = 2pt, label=170:{$15$}] (B3) at (180: 3.5) {};
\node[vtx, fill=black, inner sep = 2pt, label=180:{$7$}] (B4) at (225: 4) {};
\node[vtx, fill=black, inner sep = 2pt, label=0:{$19$}] (B5) at (315: 4) {};
\draw[thick] (B2) -- (A3) -- (B3) -- (B4) -- (A4) -- (A5) -- (A1) -- (A2) -- (B2);
\draw[thick] (A0) -- (B5) -- (B0) -- (B1) -- (A0);
\draw[dashed, thick] (A1) -- (B1) -- (B2) -- (B3);
\draw[dashed, thick] (A3) -- (B4) -- (B5) -- (A5);
\draw[dashed, thick] (A0) -- (B0);
\draw[dashed, thick] (A2) -- (A4);
\begin{scope}[on background layer]
\draw[dashed, thick, shorten <= -5mm] (B0.0) -- (B0.180);
\draw[dashed, thick, shorten <= -5mm] (A0.180) -- (A0.0);
\draw[dashed, thick, shorten <= -5mm] (B3.180) -- (B3.0);
\draw[dashed, thick, shorten <= -5mm] (A3.0) -- (A3.180);
\draw[dashed, thick, shorten <= -5mm] (A1.225) -- (A1.45);
\draw[dashed, thick, shorten <= -5mm] (A2.315) -- (A2.135);
\draw[dashed, thick, shorten <= -5mm] (A4.45) -- (A4.225);
\draw[dashed, thick, shorten <= -5mm] (A5.135) -- (A5.315);
\draw[thick, opacity = .5, blue, line width = 3pt] (A4) -- (A5);
\draw[thick, shorten <= -5mm, opacity = .5, blue, line width = 3 pt] (A4.45) -- (A4.225);
\draw[thick, shorten <= -5mm, opacity = .5, blue, line width = 3 pt] (A5.135) -- (A5.315);
\end{scope}
\end{tikzpicture}
}
\hspace{0mm}
\subfigure{
\begin{tikzpicture}[scale = 0.55]
\node[vtx, fill=black, inner sep = 2pt, label=90:{$0$}] (C) at (90: 3.5) {};
\node[vtx, fill=black, inner sep = 2pt, label=180:{$18$}] (A0) at (0: 1.9) {};
\node[vtx, fill=black, inner sep = 2pt, label=260:{$20$}] (A1) at (60: 1.9) {};
\node[vtx, fill=black, inner sep = 2pt, label=280:{$4$}] (A2) at (120: 1.9) {};
\node[vtx, fill=black, inner sep = 2pt, label=0:{$22$}] (A3) at (180: 1.9) {};
\node[vtx, fill=black, inner sep = 2pt, label=270:{$8$}] (A4) at (230: 1.9) {};
\node[vtx, fill=black, inner sep = 2pt, label=270:{$10$}] (A5) at (310: 1.9) {};
\node[vtx, fill=black, inner sep = 2pt, label=0:{$14$}] (B0) at (0: 3.5) {};
\node[vtx, fill=black, inner sep = 2pt, label=0:{$24$}] (B1) at (60: 3.3) {};
\node[vtx, fill=black, inner sep = 2pt, label=180:{$16$}] (B2) at (120: 3.3) {};
\node[vtx, fill=black, inner sep = 2pt, label=180:{$2$}] (B3) at (180: 3.5) {};
\node[vtx, fill=black, inner sep = 2pt, label=180:{$12$}] (B4) at (210: 3.3) {};
\node[vtx, fill=black, inner sep = 2pt, label=0:{$6$}] (B5) at (330: 3.3) {};
\draw[thick] (B2) -- (A1) -- (A0) -- (A5) -- (B5) -- (B0) -- (B1) -- (C) -- (B2) -- (B3) -- (A2) -- (B1);
\draw[thick] (B4) -- (A3) -- (A4);
\draw[dashed, thick] (A4) -- (B3) -- (B4) -- (A4);
\draw[dashed, thick] (A2) -- (A3) -- (B2);
\draw[dashed, thick] (A5) -- (B0) -- (A1);
\draw[dashed, thick] (B5) -- (A0) -- (B1);
\begin{scope}[on background layer]
\draw[dashed, thick, shorten <= -5mm] (B5.290) -- (B5.110);
\draw[dashed, thick, shorten <= -5mm] (B4.250) -- (B4.70);
\draw[dashed, thick, shorten <= -5mm] (A5.220) -- (A5.40);
\draw[dashed, thick, shorten <= -5mm] (A4.320) -- (A4.140);
\draw[dashed, thick, shorten <= -5mm] (A1.270) -- (A1.90);
\draw[dashed, thick, shorten <= -5mm] (A2.270) -- (A2.90);
\draw[thick, opacity = .5, blue, line width = 3pt] (A5) -- (B5);
\draw[thick, shorten <= -5mm, opacity = .5, blue, line width = 3 pt] (A5.220) -- (A5.40);
\draw[thick, shorten <= -5mm, opacity = .5, blue, line width = 3 pt] (B5.290) -- (B5.110);
\end{scope}
\end{tikzpicture}
}
\caption{Examples of orders $23$, $24$ and $25$, admitting ``infinite'' extension.}
\label{fig:n23_24_25}
\end{center}
\end{figure}

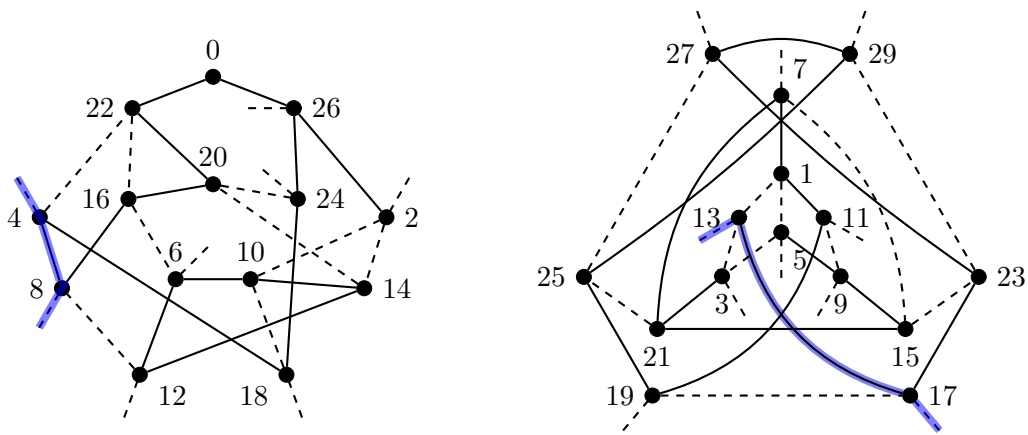
\begin{figure}[htbp]
\begin{center}
\subfigure{
\begin{tikzpicture}[scale = 0.65]
\node[vtx, fill=black, inner sep = 2pt, label=90:{$0$}] (C) at (90: 3.5) {};
\node[vtx, fill=black, inner sep = 2pt, label=180:{$22$}] (U0) at (120: 3.3) {};
\node[vtx, fill=black, inner sep = 2pt, label=0:{$26$}] (U1) at (60: 3.3) {};
\node[vtx, fill=black, inner sep = 2pt, label=180:{$4$}] (M0) at (170: 3.6) {};
\node[vtx, fill=black, inner sep = 2pt, label=180:{$16$}] (M1) at (150: 2) {};
\node[vtx, fill=black, inner sep = 2pt, label=90:{$20$}] (M2) at (90: 1.3) {};
\node[vtx, fill=black, inner sep = 2pt, label=0:{$24$}] (M3) at (30: 2) {};
\node[vtx, fill=black, inner sep = 2pt, label=0:{$2$}] (M4) at (10: 3.6) {};
\node[vtx, fill=black, inner sep = 2pt, label=180:{$8$}] (N0) at (195: 3.2) {};
\node[vtx, fill=black, inner sep = 2pt, label=90:{$6$}] (N1) at (220: 1) {};
\node[vtx, fill=black, inner sep = 2pt, label=90:{$10$}] (N2) at (320: 1) {};
\node[vtx, fill=black, inner sep = 2pt, label=0:{$14$}] (N3) at (345: 3.2) {};
\node[vtx, fill=black, inner sep = 2pt, label=350:{$12$}] (L0) at (240: 3) {};
\node[vtx, fill=black, inner sep = 2pt, label=190:{$18$}] (L1) at (300: 3) {};
\draw[thick] (M4) -- (U1) -- (M3) -- (L1) -- (M0) -- (N0) -- (M1) -- (M2) -- (U0) -- (C) -- (U1);
\draw[thick] (L0) -- (N1) -- (N2) -- (N3) -- (L0);
\draw[dashed, thick] (L0) -- (N0);
\draw[dashed, thick] (M0) -- (U0) -- (M1) -- (N1);
\draw[dashed, thick] (L1) -- (N2) -- (M4) -- (N3) -- (M2) -- (M3);
\draw[dashed, thick, shorten <= -5mm] (U1.180) -- (U1.0);
\draw[dashed, thick, shorten <= -5mm] (M4.60) -- (M4.240);
\draw[dashed, thick, shorten <= -5mm] (M3.140) -- (M3.320);
\draw[dashed, thick, shorten <= -5mm] (M0.120) -- (M0.300);
\draw[dashed, thick, shorten <= -5mm] (N0.240) -- (N0.60);
\draw[dashed, thick, shorten <= -5mm] (N1.45) -- (N1.225);
\draw[dashed, thick, shorten <= -5mm] (L0.250) -- (L0.70);
\draw[dashed, thick, shorten <= -5mm] (L1.290) -- (L1.110);
\draw[thick, opacity = .5, blue, line width = 3pt] (M0) -- (N0);
\draw[thick, shorten <= -5mm, opacity = .5, blue, line width = 3 pt] (M0.120) -- (M0.300);
\draw[thick, shorten <= -5mm, opacity = .5, blue, line width = 3 pt] (N0.240) -- (N0.60);
\begin{scope}[on background layer]
\end{scope}
\end{tikzpicture}
}
\hspace{8mm}
\subfigure{
\begin{tikzpicture}[scale = 0.65]
\node[vtx, fill=black, inner sep = 2pt, label=180:{$27$}] (O0) at (110: 4.1) {};
\node[vtx, fill=black, inner sep = 2pt, label=180:{$25$}] (O1) at (190: 4.1) {};
\node[vtx, fill=black, inner sep = 2pt, label=180:{$19$}] (O2) at (230: 4.1) {};
\node[vtx, fill=black, inner sep = 2pt, label=0:{$17$}] (O3) at (310: 4.1) {};
\node[vtx, fill=black, inner sep = 2pt, label=0:{$23$}] (O4) at (350: 4.1) {};
\node[vtx, fill=black, inner sep = 2pt, label=0:{$29$}] (O5) at (70: 4.1) {};
\node[vtx, fill=black, inner sep = 2pt, label=80:{$7$}] (M0) at (90: 3) {};
\node[vtx, fill=black, inner sep = 2pt, label=270:{$21$}] (M1) at (215: 3.1) {};
\node[vtx, fill=black, inner sep = 2pt, label=270:{$15$}] (M2) at (325: 3.1) {};
\node[vtx, fill=black, inner sep = 2pt, label=0:{$1$}] (I0) at (90: 1.4) {};
\node[vtx, fill=black, inner sep = 2pt, label=270:{$3$}] (I2) at (210: 1.4) {};
\node[vtx, fill=black, inner sep = 2pt, label=270:{$9$}] (I4) at (330: 1.4) {};
\node[vtx, fill=black, inner sep = 2pt, label=180:{$13$}] (I1) at (150: 1) {};
\node[vtx, fill=black, inner sep = 2pt, label=0:{$11$}] (I3) at (30: 1) {};
\node[vtx, fill=black, inner sep = 2pt, label=280:{$5$}] (C) at (90: 0.2) {};
\draw[thick] (I2) -- (M1) -- (M2) -- (I4) -- (C);
\draw[thick] (M0) -- (I0) -- (I3);
\draw[thick] (O2) -- (O1);
\draw[thick] (O0) to [bend left = 20] (O5);
\draw[thick] (O4) -- (O3);
\draw[thick] (O3) to [bend left = 30] (I1);
\draw[thick] (I3) to [bend left = 30] (O2);
\draw[thick] (O5) to [bend left = 5] (O1);
\draw[thick] (O4) to [bend left = 5] (O0);
\draw[thick] (M1) to [bend left = 25] (M0);
\draw[dashed, thick] (C) -- (I2) -- (I1) -- (I0) -- (C);
\draw[dashed, thick] (I4) -- (I3);
\draw[dashed, thick] (O2) -- (O3);
\draw[dashed, thick] (M2) -- (O4) -- (O5);
\draw[dashed, thick] (M1) -- (O1) -- (O0);
\draw[dashed, thick] (M0) to [bend left = 25] (M2);
\begin{scope}[on background layer]
\draw[dashed, thick, shorten <= -5mm] (C.270) -- (C.90);
\draw[dashed, thick, shorten <= -5mm] (I2.300) -- (I2.120);
\draw[dashed, thick, shorten <= -5mm] (I4.240) -- (I4.60);
\draw[dashed, thick, shorten <= -5mm] (I1.210) -- (I1.30);
\draw[dashed, thick, shorten <= -5mm] (I3.330) -- (I3.150);
\draw[dashed, thick, shorten <= -5mm] (M0.90) -- (M0.270);
\draw[dashed, thick, shorten <= -5mm] (O2.230) -- (O2.50);
\draw[dashed, thick, shorten <= -5mm] (O3.310) -- (O3.130);
\draw[dashed, thick, shorten <= -5mm] (O0.110) -- (O0.290);
\draw[dashed, thick, shorten <= -5mm] (O5.70) -- (O5.250);
\draw[thick, opacity = .5, blue, line width = 3pt] (O3) to [bend left = 30] (I1);
\draw[thick, shorten <= -5mm, opacity = .5, blue, line width = 3 pt] (O3.310) -- (O3.130);
\draw[thick, shorten <= -5mm, opacity = .5, blue, line width = 3 pt] (I1.210) -- (I1.30);
\end{scope}
\end{tikzpicture}
}
\caption{Examples of orders $27$ and $30$, admitting ``infinite'' extension.}
\label{fig:n27_30}
\end{center}
\end{figure}

\begin{proof}
In view of Proposition~\ref{pro:wreath}, the claim from (i) about even orders is clear. We confirm that it also holds for all odd orders $n$ with $n \geq 21$ by finding one ``nice enough'' example for each of the orders $n \in \{21, 23, 25, 27\}$, which allows us to apply Corollary~\ref{cor:SRcon} using $\Wr(4)$ as $\G'$ as many times as desired. To explain what exactly we mean by this consider again the self-reverse distance magic labeling of $\Wr(4)$ from Figure~\ref{fig:W(4)quo} and note that the corresponding quotient has four solid edges (whose endvertices all have semiedges). For two opposite edges in the quotient, the labels of their endvertices differ by $2$, while for the other two edges they differ by $4$. Now, suppose we have a self-reverse distance magic labeling of some tetravalent graph $\G$ of order $n$ such that in the corresponding quotient we have a solid edge with semiedges at its endvertices and such that their labels differ by $4$. We can then set $\G' = \Wr(4)$, use the above mentioned labeling of $\Wr(4)$ and take the solid edge in the quotient whose endvertices have labels $1$ and $5$. Applying Corollary~\ref{cor:SRcon}, we obtain a graph of order $8+n$ and a corresponding self-reverse distance magic labeling. Moreover, in the corresponding quotient, the vertices corresponding to the vertices labeled $3$ and $7$ in the quotient of $\G'$ are still connected by a solid edge, have semiedges, and have labels $3+n$ and $7+n$. We can thus repeat the whole process using $\G' = \Wr(4)$ and Corollary~\ref{cor:SRcon}. In this way we get one example for each of the orders of the form $n + 8m$, $m \geq 0$. Observe that if $\G$ is connected, then so is the resulting graph $\G\oplus^{C}_{C'} \Wr(4)$. Moreover, by Proposition~\ref{pro:degenerate}, the resulting distance magic labeling is not degenerate (since the vertex labeled $n+7$ is adjacent to the one labeled $n+3$, but not to the one labeled $-n-3$).

Since by Proposition~\ref{pro:degenerate}, a self-reverse distance magic labeling of a connected tetravalent graph of odd order (which is of course not a wreath graph) is automatically non-degenerate, this will also prove parts (ii) and (iii) for odd orders $n$. By Proposition~\ref{pro:wreath} and the data from Table~\ref{tab:data}, the parts (ii) and (iii) for $n$ even will be confirmed if we can provide one example for each of the orders $n \in \{18, 20, 24, 30\}$ such that the corresponding graph is not a wreath graph and the corresponding quotient has a solid edge allowing an extension using Corollary~\ref{cor:SRcon} as in the previous paragraph. The examples (given in terms of the corresponding quotient with an ``appropriate'' solid edge and the corresponding semiedges highlighted in blue) for each of $n \in \{18,20,21,23,24,25,27,30\}$ are presented in Figures~\ref{fig:n18_20_21},~\ref{fig:n23_24_25}, and~\ref{fig:n27_30}.
\end{proof}


\section{Examples with a large degree of symmetry and future directions}
\label{sec:VT}

In this last section of the paper we focus on tetravalent distance magic graphs with a large degree of symmetry. One might argue that this is not interesting since the properties of a graph being distance magic and it being vertex-transitive (see the next paragraph for a definition) seem to be completely unrelated. However, the scarcity of examples of vertex-transitive distance magic graphs (see the rest of this section) makes such an investigation very appealing. Moreover, let us mention that the well-known open question of Lov\'asz from 1969~\cite{Lov69}, asking whether each connected vertex-transitive graph possesses a Hamilton path, also links two seemingly unrelated concepts in graph theory. Nevertheless, the question has spurred a tremendous amount of interest and research (see for instance the survey~\cite{KutMar09} or~\cite{KutMarRaz25} for a list of some more recent references). We think that also in our case it is interesting to investigate in what way, if at all, the property of being distance magic is related to the degree of symmetry of a graph.

Now, recall first that a graph is {\em vertex-transitive} or {\em edge-transitive}, if its automorphism group acts transitively on the vertex set or the edge set of the graph, respectively. It is well known and easy to see that the wreath graphs $\Wr(m)$ are vertex-transitive and also edge-transitive. We were thus curious ``how symmetric'' the connected tetravalent graphs admitting a non-degenerate self-reverse distance magic labeling are. The last line of the above Table~\ref{tab:data} contains the number of connected tetravalent vertex-transitive graphs of a given order that admit a non-degenerate self-reverse distance magic labeling ($\# \text{VT}$). The obtained data suggests that there seem to be very few examples of such graphs. In particular, no such graph of odd order up to $29$ exists. In some sense, this may not seem very surprising as one might think that the central vertex should be ``special''. However, the central vertex is only central with respect to a given labeling. Moreover, there exist many tetravalent vertex-transitive (but not necessarily distance magic) graphs of odd order with the property that for each vertex an involutory automorphism fixing precisely this vertex (and no other) exists --- each circulant of odd order has this property (see the last paragraph of this section for a definition of a circulant). It is thus not clear whether connected tetravalent vertex-transitive graphs of odd order admitting a self-reverse distance magic labeling might exist. In fact, we are not aware of any known example of a connected tetravalent distance magic vertex-transitive graph of odd order. This thus suggests the following question.

\begin{question}
\label{que:VTodd}
Do there exist connected tetravalent vertex-transitive distance magic graphs of odd order? If so, does any of them admit a self-reverse distance magic labeling?
\end{question}

The situation for even order is of course quite different, since for each integer $m \geq 3$, the wreath graph $\Wr(m)$ admits a self-reverse distance magic labeling. However, if we exclude the wreath graphs from our consideration, then the situation again becomes very interesting. Namely, the data from Table~\ref{tab:data} reveals that up to order $30$, only five connected tetravalent ``non-wreath'' vertex-transitive graphs admitting a self-reverse distance magic labeling exist (recall that by Proposition~\ref{pro:wreath}, the vertex-transitive example of order $16$ from Table~\ref{tab:data} is $\Wr(8)$ and one of the examples of order $24$ is $\Wr(12)$). Two are of order $24$, while we have one for each of the orders $18$, $20$ and $30$.

\begin{figure}[htbp]
\begin{center}
\subfigure{
\begin{tikzpicture}[scale = 0.65]
\node[vtx, fill=black, inner sep = 2pt, label=80:{$1$}] (O0) at (90: 4) {};
\node[vtx, fill=black, inner sep = 2pt, label=190:{$7$}] (O1) at (162: 4) {};
\node[vtx, fill=black, inner sep = 2pt, label=280:{$11$}] (O2) at (234: 4) {};
\node[vtx, fill=black, inner sep = 2pt, label=260:{$5$}] (O3) at (306: 4) {};
\node[vtx, fill=black, inner sep = 2pt, label=350:{$9$}] (O4) at (18: 4) {};
\node[vtx, fill=black, inner sep = 2pt, label=0:{$17$}] (I0) at (90: 2.3) {};
\node[vtx, fill=black, inner sep = 2pt, label=270:{$19$}] (I1) at (162: 2.3) {};
\node[vtx, fill=black, inner sep = 2pt, label=280:{$13$}] (I2) at (234: 2.3) {};
\node[vtx, fill=black, inner sep = 2pt, label=260:{$3$}] (I3) at (306: 2.3) {};
\node[vtx, fill=black, inner sep = 2pt, label=270:{$15$}] (I4) at (18: 2.3) {};
\draw[thick] (O0) -- (I0) -- (I2) -- (O2) -- (O3) -- (I3) -- (I0);
\draw[thick] (O1) -- (I1) -- (I4) -- (O4);
\draw[dashed, thick] (O2) -- (O1) -- (O0) -- (O4) -- (O3);
\draw[dashed, thick] (I1) -- (I3);
\draw[dashed, thick] (I2) -- (I4);
\begin{scope}[on background layer]
\draw[dashed, thick, shorten <= -5mm] (O0.90) -- (O0.270);
\draw[dashed, thick, shorten <= -5mm] (O1.162) -- (O1.342);
\draw[dashed, thick, shorten <= -5mm] (O2.234) -- (O2.54);
\draw[dashed, thick, shorten <= -5mm] (O3.306) -- (O3.126);
\draw[dashed, thick, shorten <= -5mm] (O4.18) -- (O4.198);
\draw[dashed, thick, shorten <= -5mm] (I0.270) -- (I0.90);
\draw[dashed, thick, shorten <= -5mm] (I1.342) -- (I1.162);
\draw[dashed, thick, shorten <= -5mm] (I2.54) -- (I2.234);
\draw[dashed, thick, shorten <= -5mm] (I3.126) -- (I3.306);
\draw[dashed, thick, shorten <= -5mm] (I4.198) -- (I4.18);

\end{scope}
\end{tikzpicture}
}
\hspace{8mm}
\subfigure{
\begin{tikzpicture}[scale = 0.65]
\node[vtx, fill=black, inner sep = 2pt, label=180:{$1$}] (O0) at (150: 4.5) {};
\node[vtx, fill=black, inner sep = 2pt, label=270:{$12$}] (O1) at (270: 4.5) {};
\node[vtx, fill=black, inner sep = 2pt, label=0:{$9$}] (O2) at (30: 4.5) {};
\node[vtx, fill=black, inner sep = 2pt, label=90:{$8$}] (M0) at (120: 2.5) {};
\node[vtx, fill=black, inner sep = 2pt, label=180:{$11$}] (M1) at (180: 2.5) {};
\node[vtx, fill=black, inner sep = 2pt, label=180:{$10$}] (M2) at (240: 2.5) {};
\node[vtx, fill=black, inner sep = 2pt, label=0:{$2$}] (M3) at (300: 2.5) {};
\node[vtx, fill=black, inner sep = 2pt, label=0:{$6$}] (M4) at (0: 2.5) {};
\node[vtx, fill=black, inner sep = 2pt, label=90:{$7$}] (M5) at (60: 2.5) {};
\node[vtx, fill=black, inner sep = 2pt, label=270:{$5$}] (I0) at (90: 1.3) {};
\node[vtx, fill=black, inner sep = 2pt, label=10:{$3$}] (I1) at (210: 1.3) {};
\node[vtx, fill=black, inner sep = 2pt, label=170:{$4$}] (I2) at (330: 1.3) {};
\draw[thick] (O0) to [bend left = 30] (O2);
\draw[thick] (O2) to [bend left = 30] (O1);
\draw[dashed, thick] (O1) to [bend left = 30] (O0);
\draw[thick] (O1) -- (M3) -- (I2) -- (M4) -- (M5) -- (M0) -- (I0) -- (I1) -- (M2) -- (M1) -- (O0);
\draw[dashed, thick] (O1) -- (M2) -- (M3) -- (M4) -- (O2) -- (M5) -- (I0) -- (I2) -- (I1) -- (M1) -- (M0) -- (O0);
\end{tikzpicture}
}
\caption{The ``curious'' vertex-transitive examples of orders $20$ and $24$.}
\label{fig:VTexamples}
\end{center}
\end{figure}
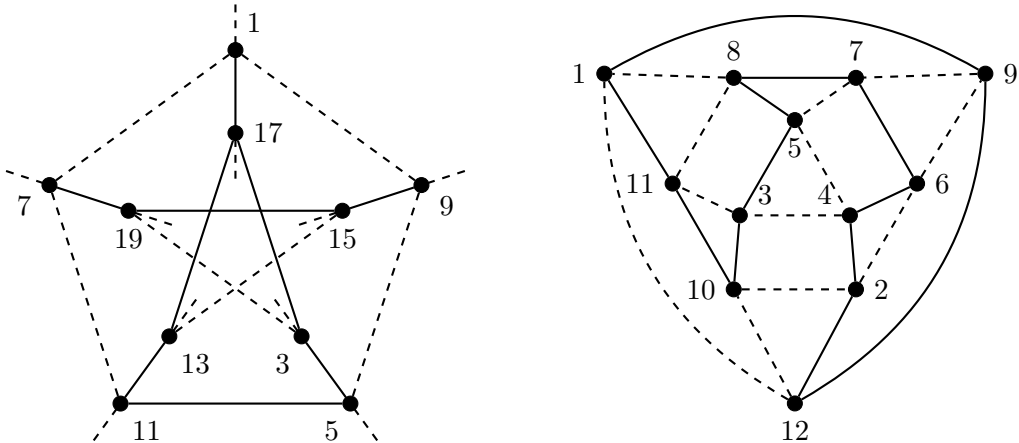

Three of these belong to well-known families of tetravalent distance magic graphs. The unique example of order $18$ is the Cartesian product $C_3 \square C_6$ (see~\cite{RozSpa24} for the classification of all distance magic Cartesian products of two cycles). One of the two examples of order $24$ is the circulant $\mathrm{Circ}(24; \{\pm 1, \pm 5\})$ and the unique one of order $30$ is the circulant $\mathrm{Circ}(30; \{\pm 1, \pm 4\})$ (see \cite{CicFro16, MikSpa21} for a complete classification of tetravalent distance magic circulants). The remaining two graphs, the unique example of order $20$ and the remaining example of order $24$, are presented in Figure~\ref{fig:VTexamples} (via the corresponding quotient for one of the possible self-reverse distance magic labelings).

It is not difficult to show (but one can simply use a computer) that of these five examples only $\mathrm{Circ}(24; \{\pm 1, \pm 5\})$ is edge-transitive. The results of~\cite{CicFro16, MikSpa21} in fact imply that each distance magic circulant of the form $\mathrm{Circ}(n ; \{\pm 1, \pm c\})$, where $c$ is odd, is edge-transitive. Moreover, by the results of~\cite{AnhCicPetTep15}, the direct product of a cycle $C_m$ with itself in the case that $m$ is divisible by $4$ is distance magic (it of course consists of two components) and it is clearly vertex- and edge-transitive. However, since we are aware of no other known tetravalent edge-transitive (or vertex-transitive) distance magic graphs, the following two problems seem very natural (the first is a generalization of~\cite[Problem~5.5]{MikSpa21}).

\begin{problem}
Characterize connected tetravalent edge-transitive distance magic graphs and determine which of them admit self-reverse distance magic labelings.
\end{problem}

\begin{problem}
Characterize connected tetravalent vertex-transitive distance magic graphs and determine which of them admit self-reverse distance magic labelings.
\end{problem}

One of the first natural steps to be taken regarding the above two problems is for sure a thorough investigation of the infinite families of tetravalent graphs whose members are known to be distance magic in order to determine which of these graphs admit self-reverse distance magic labelings. We already have some results in this regard but since the main purpose of this paper was to introduce the concept of self-reverse distance magic labelings, we do not want to make it any longer by presenting the corresponding proofs. Let us simply say that the results of~\cite{RozSpa24} imply that each distance magic labeling of a Cartesian product of two cycles is in fact self-reverse and that we can prove that each tetravalent distance magic circulant admits a self-reverse distance magic labeling. We also have a proof that for each distance magic direct product of two cycles (see~\cite{AnhCicPetTep15} for the classification) a connected component of this graph (which is of course vertex-transitive) admits a self-reverse distance magic labeling. This thus suggests another natural question.

\begin{question}
Does there exist a connected tetravalent vertex-transitive distance magic graph that admits no self-reverse distance magic labeling?
\end{question}

As concerns a search for possible new families of tetravalent vertex-transitive (or even edge-transitive) distance magic graphs, the extensive literature on the theory of vertex- and/or edge-transitive tetravalent graphs offers numerous possibilities. It was suggested already in~\cite{MikSpa21} that the well-known Rose window graphs~\cite{Wil08} seem like a very good family to start with. 

Let us conclude the paper with yet another interesting aspect of the obtained data presented in Table~\ref{tab:data}. The wreath graphs, the circulants and the Cartesian and direct products of two cycles are of course all Cayley graphs. For self-completeness we recall that for a group $G$ and its inverse-closed subset $S$ not containing the identity the {\em Cayley graph} of $G$ {\em with respect} to $S$ is the graph with vertex set $G$ in which vertices $g$ and $h$ are adjacent whenever $g^{-1}h \in S$. In the case that $G$ is cyclic the graph is often called a {\em circulant}. Now, on top of what was said about the known families of tetravalent distance magic graphs, we point out that it can be verified that the example of order $24$, presented on the right part of Figure~\ref{fig:VTexamples}, is a Cayley graph (of the symmetric group $S_4$). On the other hand, it should not be surprising that the example presented on the left part of Figure~\ref{fig:VTexamples} is not a Cayley graph (just like the Petersen graph is not), and so this is to the best of our knowledge the only known connected tetravalent vertex-transitive graph which is not a Cayley graph but admits a self-reverse distance magic labeling. In fact, we know of no other connected tetravalent vertex-transitive distance magic graph which is not a Cayley graph. This brings us to the following natural question.

\begin{question}
\label{que:VTnCay}
Are there infinitely many connected tetravalent vertex-transitive distance magic graphs which are not Cayley graphs? If so, are there infinitely many of them that admit a self-reverse distance magic labeling?
\end{question}

\section*{Declarations and Acknowledgments}

\noindent
{\bf Competing Interests:} The authors have no relevant competing interests to declare.
\bigskip

\noindent
{\bf Funding:} Petr Kov\'{a}\v{r} was co-funded by the financial support of the European Union under the REFRESH -- Research Excellence For Region Sustainability and High-tech Industries project number CZ.10.03.01/00/22 003/0000048 via the Operational Programme Just Transition. Ksenija Rozman and Primo\v z \v Sparl acknowledge the financial support by the Slovenian Research and Innovation Agency (research program P1-0285 and research projects J1-3001 and J1-50000).

\end{document}